# STATISTICAL INFERENCE FOR SEMIPARAMETRIC VARYING-COEFFICIENT PARTIALLY LINEAR MODELS WITH ERROR-PRONE LINEAR COVARIATES

BY YONG ZHOU[1] AND HUA LIANG[2]

*Shanghai University of Finance and Economics, Chinese Academy of Sciences and University of Rochester*

We study semiparametric varying-coefficient partially linear models when some linear covariates are not observed, but ancillary variables are available. Semiparametric profile least-square based estimation procedures are developed for parametric and nonparametric components after we calibrate the error-prone covariates. Asymptotic properties of the proposed estimators are established. We also propose the profile least-square based ratio test and Wald test to identify significant parametric and nonparametric components. To improve accuracy of the proposed tests for small or moderate sample sizes, a wild bootstrap version is also proposed to calculate the critical values. Intensive simulation experiments are conducted to illustrate the proposed approaches.

**1. Introduction.** Various efforts have been made to balance the interpretation of linear models and flexibility of nonparametric models. Important results from these efforts include semiparametric varying-coefficient partially linear models (SVCPLM), in which the response variable $Y$ depends on variables $\mathbf{Z}$, $\mathbf{X}$ and $U$ in the form of

$$(1.1) \qquad Y = \mathbf{\Theta}^{\mathrm{T}}\mathbf{Z} + \boldsymbol{\alpha}^{\mathrm{T}}(U)\mathbf{X} + \varepsilon,$$

Received July 2007.
[1]Supported in part by National Natural Science Foundation of China (NSFC) Grant (No. 10471140 and No. 10731010), the National Basic Research Program (No. 2007CB814902) and Creative Research Groups of China, SUFE funding through Project 211 Phase III and Shanghai Leading Academic Discipline Project (B803).
[2]Supported by NIH/NIAID Grants AI62247 and AI59773 and NSF grant DMS806097.

where $\boldsymbol{\Theta}$ is a $p$-dimensional vector of unknown parameters, $\boldsymbol{\alpha}(\cdot)$ is a $q$-variate vector of unknown functions, $U$ is a vector of nonparametric components that may be multivariate and the model error $\varepsilon$ has mean zero and finite variance. For notational simplicity, we assume that $U$ is scalar. $\boldsymbol{\alpha}^{\mathrm{T}}(U)\mathbf{X}$ is referred to as a nonparametric component since $\boldsymbol{\alpha}(U)$ is nonparametric.

Model (1.1) permits the interaction between the covariates $U$ and $\mathbf{X}$ in such a way that a different level of covariate $U$ is associated with a different linear model about $\boldsymbol{\Theta}^{\mathrm{T}}\mathbf{Z}$, and allows one to examine the extent to which covariates $\mathbf{X}$ interact. This model presents a novel and general structure, which indeed covers many well-studied, important semiparametric regression models. For example, when $\mathbf{Z} = 0$, (1.1) reduces to varying-coefficient models, which were originally proposed by Hastie and Tibshirani (1993) and studied by Fan and Zhang (1999), Xia and Li (1999) and Cai, Fan and Li (2000). When $q = 1$ and $\mathbf{X} = 1$, (1.1) reduces to well-known partially linear models, in which $Y$ depends on $\mathbf{Z}$ in a linear way but is related to another independent variable $U$ in an unspecified form. There is a great deal of literature on the study of partially linear models [e.g., Engle et al. (1986), Robinson (1988) and Speckman (1988)]. A survey of partially linear models was given by Härdle, Liang and Gao (2000). The study of SVCPLM has been investigated by Zhang, Lee and Song (2002) and Fan and Huang (2005), among others. Zhang, Lee and Song (2002) developed the procedures for estimation of the linear and nonparametric parts of the SVCPLM. Fan and Huang (2005) proposed a profile likelihood technique for estimating parametric components and established the asymptotic normality of their proposed estimator.

All studies of the SVCPLM are limited to considerations of exactly observed data. However, in biomedical research observations are measured with error. Simply ignoring measurement errors, known as the naive method, will result in biased estimators. Various attempts have been made to correct for such bias, see Fuller (1987) and Carroll et al. (2006) for extensive discussions and examples of linear and nonlinear models with measurement errors. In this paper, we are concerned with the situation where some components ($\boldsymbol{\xi}$) of $\mathbf{Z}$ are unobserved directly, but auxiliary information is available to remit $\boldsymbol{\xi}$. Let $\mathbf{Z} = (\boldsymbol{\xi}^{\mathrm{T}}, \mathbf{W}^{\mathrm{T}})^{\mathrm{T}}$, where $\boldsymbol{\xi}$ is a $p_1 \times 1$ vector and $\mathbf{W}$ is a vector of the remaining observed components. We assume that $\boldsymbol{\xi}$ is related to observed $\boldsymbol{\eta}$ and $\mathbf{V}$ through the relationship $\boldsymbol{\xi} = E(\boldsymbol{\eta}|\mathbf{V})$. Thus, we study the following model:

$$(1.2) \qquad \begin{cases} Y = \boldsymbol{\beta}^{\mathrm{T}}\boldsymbol{\xi} + \boldsymbol{\theta}^{\mathrm{T}}\mathbf{W} + \boldsymbol{\alpha}^{\mathrm{T}}(U)\mathbf{X} + \varepsilon, \\ \boldsymbol{\eta} = \boldsymbol{\xi}(\mathbf{V}) + \mathbf{e}, \end{cases}$$

where $E(\varepsilon|\mathbf{Z}, \mathbf{X}, U) = 0$, $E(\varepsilon^2|\mathbf{Z}, \mathbf{X}, U) = \sigma^2(\mathbf{Z}, \mathbf{X}, U)$ and $\mathbf{e}$ is an error with mean zero and positive finite covariance matrix $\Sigma_e = E(\mathbf{e}\mathbf{e}^{\mathrm{T}})$. The four covariates $\mathbf{V}$, $\mathbf{W}$, $\mathbf{X}$ and $U$ are different. In our structure, we allow that $\mathbf{V}$ and

$(\mathbf{X}, \mathbf{W}, U)$ may overlap. Model (1.2) is flexible enough to include a variety of models of interest. We give three examples to illustrate its flexibility:

EXAMPLE 1 (Errors-in-variable models with validation data). $\mathbf{Z}$ is a $p$-variate variable vector and is not observed. $\tilde{\mathbf{Z}}$ is an another $p$-variate vector and is observed associated with vector $\mathbf{Z}$. Assume that we have primary observations $\{Y_j, \tilde{\mathbf{Z}}_j, U_j, j = 1, \ldots, n\}$, and $n_0$ independent validation observations $\{\mathbf{Z}_j, \tilde{\mathbf{Z}}_j, U_i, i = n+1, \ldots, n+n_0\}$, which are independent of the primary observations. Let $\mathbf{V} = (\tilde{\mathbf{Z}}^{\mathrm{T}}, U)^{\mathrm{T}}$. The partial errors-in-variable model with validation data is written as

$$(1.3) \qquad \begin{cases} Y = \boldsymbol{\beta}^{\mathrm{T}} E(\mathbf{Z}|\mathbf{V}) + \alpha(U) + \varepsilon, \\ \varepsilon = e + \boldsymbol{\beta}^{\mathrm{T}}\{\mathbf{Z} - E(\mathbf{Z}|\mathbf{V})\}. \end{cases}$$

This model has been studied by Sepanski and Lee (1995), Sepanski and Carroll (1993) and Sepanski, Knickerbocker and Carroll (1994). Taking $\mathbf{X} = 1$, $\boldsymbol{\theta} = 0$, $\eta = \mathbf{Z}$ and $\boldsymbol{\xi} = E(\mathbf{Z}|\mathbf{V})$ in (1.2), we know that (1.3) is a sub-model of (1.2).

EXAMPLE 2 (De-noise linear model). The relation between the response variable $Y$ and covariates $(\boldsymbol{\xi}, \mathbf{W})$ is described by $Y = \boldsymbol{\beta}^{\mathrm{T}}\boldsymbol{\xi} + \boldsymbol{\theta}^{\mathrm{T}}\mathbf{W} + \varepsilon$, where $\boldsymbol{\beta}$ and $\boldsymbol{\theta}$ are parametric vectors, respectively. The covariate $\boldsymbol{\xi}$ is measured with error since, instead of observing $\boldsymbol{\xi}$ directly, we observe its surrogate $\boldsymbol{\eta}$. This forms a de-noise linear model:

$$(1.4) \qquad \begin{cases} Y = \boldsymbol{\beta}^{\mathrm{T}}\boldsymbol{\xi} + \boldsymbol{\theta}^{\mathrm{T}}\mathbf{W} + \varepsilon, \\ \boldsymbol{\eta} = \boldsymbol{\xi} + e, \end{cases}$$

where $\boldsymbol{\xi} = \boldsymbol{\xi}(t)$ is subject to measurement error at time $t$ and the measurement errors $\varepsilon$ and $e$ are independent of each other at each time $t$.

Cai, Naik and Tsai (2000) used this model to estimate the relationship between awareness and television rating points of TV commercials for certain products. Cui, He and Zhu (2002) proposed an estimator of the coefficients and established asymptotic results of the proposed estimator. It is easy to see that (1.2) includes (1.4).

EXAMPLE 3 (Rational expectation model). Consider the following rational expectation model:

$$(1.5) \qquad Y_t = \boldsymbol{\gamma}^{\mathrm{T}}\mathbf{S}_t + \boldsymbol{\zeta}^{\mathrm{T}}\{\boldsymbol{\eta}_t - E(\boldsymbol{\eta}_t|\mathbf{V}_t)\} + \varepsilon_t,$$

where $\boldsymbol{\eta}_t - E(\boldsymbol{\eta}_t|\mathbf{V}_t)$ is the expectation payoff for price variable $\boldsymbol{\eta}_t$ given historical information $V_t$. In this model, $(Y_t, \mathbf{S}_t, \boldsymbol{\eta}_t, V_t)$ except $E(\boldsymbol{\eta}_t|\mathbf{V}_t)$ can be observed directly.

Besides estimation and inference of $\gamma$ and $\boldsymbol{\zeta}$, within the econometric community, the following model is of interest:

$$(1.6) \qquad Y_t = \gamma^{\mathrm{T}} \mathbf{S}_t + \boldsymbol{\zeta}^{\mathrm{T}} \boldsymbol{\eta}_t - \boldsymbol{\beta}^{\mathrm{T}} E(\boldsymbol{\eta}_t | \mathbf{V}_t) + \varepsilon_t.$$

It is worthy to note that (1.6) is a sub-model of (1.2). An interesting question is to test whether the (1.6) satisfies the rational expectation model (1.5), that is, to test following hypothesis:

$$(1.7) \qquad H_0 : \boldsymbol{\beta} = \boldsymbol{\zeta} \quad \mathrm{VS} \quad H_1 : \boldsymbol{\beta} \neq \boldsymbol{\zeta}.$$

In the econometric literature, the regression of unobserved covariates is also called generated regression. This topic has been widely studied. Pagan (1984) gave a comprehensive review on the estimation of parametric models with generated regression. Ai and Mcfadden (1997) presented a procedure for analyzing a partially specified nonlinear regression model in which the nuisance parameter is an unrestricted function of a subset of regressors. Ahn and Powell (1993) and Powell (1987) considered the case with the generated regressors in the nonparametric part of the model. Li (2002) considered the problems of estimating a semiparametric partially linear model for dependent data with generated regressors. Their models are special cases of the rational expectation model.

Various procedures similar to generated regression have been proposed to reduce the bias due to mismeasurement. Regression calibration and simulation extrapolation have been developed for measurement errors models Carroll et al. (2006). Liang, Härdle and Carroll (1999) studied a special case of (1.2), partially linear errors-in-variables models, and proposed an attenuated estimator of the parameter based on the semiparametric likelihood estimate. Wang and Pepe (2000) used a pseudo-expected estimating equation method to estimate the parameter in order to correct the estimation bias.

In an attempt to develop a unified estimation procedure for (1.2), we propose a profile-based procedure, which is similar to regression calibration method in spirit. The procedure consists of two steps. In the first step, we calibrate the error-prone covariate $\boldsymbol{\xi}$ by using ancillary information and applying nonparametric regression techniques. In the second step, we use profile least-square-based principle for estimating the parametric and nonparametric components. Under the mild assumptions, we derive the asymptotic representives of the proposed estimators, and use the representives to establish asymptotic normality. We also propose the profile least-square-based ratio test and Wald test for the parametric part of (1.2), and a goodness-of-fit test for the varying coefficients in the nonparametric part. The asymptotic distribution of the proposed test statistics are derived. Wild bootstrap versions are introduced to calculate the critical values for those tests.

The paper is organized as follows: In Section 2, we focus on the estimation of the parameters and nonparametric functions, and on the development of asymptotic properties of the resulting estimators. The error-prone covariates are first calibrated. Bandwidth selection strategy is also discussed. In Section 3, we develop profile least-square-based ratio tests for parametric and nonparametric components. Wild bootstrap methods are proposed to calculate the critical values. The results of applications to simulated and real data are reported in Section 4. Section 5 gives a conclusion. Regularity assumptions and technical proofs are relegated to the Appendix.

**2. Estimation of the parametric and nonparametric components.** When $\boldsymbol{\xi}$ is observed, estimators of $\boldsymbol{\beta}$ and $\alpha(u)$ and associated tests have been developed to study (1.2). These estimators and tests cannot be used directly when $\boldsymbol{\xi}$ is unobservable. We first need to calibrate $\boldsymbol{\xi}$ by using ancillary variables $\boldsymbol{\eta}$ and $\mathbf{V}$ because a direct replacement of $\boldsymbol{\xi}$ by $\boldsymbol{\eta}$ will result in bias.

2.1. *Covariate calibration.* For notational simplicity, we assume $\mathbf{V}$ is univariate in the remainder of this paper. Let $\eta_{i,k}$ be the $k$th entry of vector $\boldsymbol{\eta}$, and $L_b(\cdot) = L(\cdot/b)/b$, $b = b_k$ ($k = 1, 2, \ldots, p_1$) is a bandwidth for the $k$th component of $\boldsymbol{\eta}$. Assume throughout the paper that $\xi_k(v)$ has $r+1$ derivatives and we approximate $\xi_k(v)$ by an $r$-order polynomial within the neighborhood of $v_0$ via Taylor expansion

$$\xi_k(v) \approx \xi_k(v_0) + \xi_k'(v_0)(v-v_0) + \cdots + \frac{\xi_k^{(r)}(v_0)}{r!}(v-v_0)^r = \sum_{j=0}^r a_{j,k}(v-v_0)^j.$$

Denote

$$\mathbf{V}_v = \begin{bmatrix} 1 & (V_1 - v) & \cdots & (V_1 - v)^r \\ \vdots & \vdots & \cdots & \cdots \\ 1 & (V_n - v) & \cdots & (V_n - v)^r \end{bmatrix}, \qquad \boldsymbol{\eta}^{(k)} = \begin{pmatrix} \eta_{1k} \\ \vdots \\ \eta_{nk} \end{pmatrix},$$

$W_v = \operatorname{diag}\{L_b(V_1 - v), \ldots, L_b(V_n - v)\}$. The local polynomial estimator [Fan and Gijble (1996)] of $(a_{0,k}, \ldots, a_{r,k})^{\mathrm{T}}$ can be expressed as $\hat{\mathbf{a}}_k^{\mathrm{T}} = (\mathbf{V}_v^{\mathrm{T}} W_v \mathbf{V}_v)^{-1} \mathbf{V}_v^{\mathrm{T}} \times W_v \boldsymbol{\eta}^{(k)}$. As a consequence, $\xi_k(v)$ is estimated by $\hat{\xi}_k(v) = \zeta_1^{\mathrm{T}} (\mathbf{V}_v^{\mathrm{T}} W_v \mathbf{V}_v)^{-1} \mathbf{V}_v^{\mathrm{T}} \times W_v \boldsymbol{\eta}^{(k)}$, for $k = 1, \ldots, p_1$, where $\zeta_1$ is a $(r+1) \times 1$ vector with 1 in the first position and 0 in other positions.

In what follows, denote $A^{\otimes 2} = AA^{\mathrm{T}}$, $\mu_j = \int u^j L(u)\,du$, $\nu_j = \int u^j L^2(u)\,du$, $S_u = (\mu_{j+l})_{0 \le j, l \le r}$ and $c_p = (\mu_{r+1}, \ldots, \mu_{2r+1})^{\mathrm{T}}$. $f_v(v)$ is the density function of $V$.

Under the assumptions given in the Appendix, we can prove [Fan and Gijbels (1996), pages 101–103 or Carroll et al. (1997), page 486] that

$$\hat{\boldsymbol{\xi}}(v) - \boldsymbol{\xi}(v) = \frac{\zeta_1 S_u^{-1} c_p b^{r+1}}{(r+1)!} \boldsymbol{\xi}^{(r+1)}(v) + \frac{1}{n f_v(v)} \sum_{i=1}^n L_b(V_i - v)\mathbf{e}_i$$

(2.1)
$$+ o(b^{r+1} + \log b^{-1}/\sqrt{nb}),$$

uniformly on $v \in \mathcal{V}$. This fact will be used for proving the main results in the Appendix.

2.2. *Estimation of the parametric component.* Let $(Y_i, \boldsymbol{\eta}_i, \mathbf{V}_i, \mathbf{W}_i, \mathbf{X}_i, U_i)$, $i = 1, 2, \ldots, n$, be the observations from (1.2). The unknown covariates $\boldsymbol{\xi}_i$ are substituted by their estimators given in the above section. We therefore have following "new" model:

(2.2)
$$\begin{cases} Y_i = \boldsymbol{\beta}^{\mathrm{T}} \hat{\boldsymbol{\xi}}_i + \boldsymbol{\theta}^{\mathrm{T}} \mathbf{W}_i + \boldsymbol{\alpha}^{\mathrm{T}}(U_i)\mathbf{X}_i + \hat{\varepsilon}_i, \\ \hat{\varepsilon}_i = \varepsilon_i + \boldsymbol{\beta}^{\mathrm{T}}\{\boldsymbol{\xi}_i - \hat{\boldsymbol{\xi}}(\mathbf{V}_i)\}, \end{cases} \quad i = 1, \ldots, n,$$

where $\{\hat{\varepsilon}_i\}_{i=1}^n$ are still treated as errors. If $\hat{\boldsymbol{\xi}}_i$ would be an unbiased estimator of $\boldsymbol{\xi}_i$, then $E\hat{\varepsilon}_i = 0$.

Approximate $\alpha_j(U)$ within the neighbors of $u$ by $a_j(u) + b_j(u)(U - u)$ for $j = 1, \ldots, q$. Write $\hat{\mathbf{Z}}_i = (\hat{\boldsymbol{\xi}}_i^{\mathrm{T}}, \mathbf{W}_i^{\mathrm{T}})^{\mathrm{T}}$ and $\boldsymbol{\Theta} = (\boldsymbol{\beta}^{\mathrm{T}}, \boldsymbol{\theta}^{\mathrm{T}})^{\mathrm{T}}$. Following the profile likelihood-based procedure proposed by Fan and Huang (2005), our profile least-square-based estimator of $\boldsymbol{\Theta}$ is defined as

(2.3)
$$\hat{\boldsymbol{\Theta}}_n = \{\widetilde{\mathbf{Z}}^{\mathrm{T}} \widetilde{\mathbf{Z}}\}^{-1} \widetilde{\mathbf{Z}}^{\mathrm{T}} (\mathbf{I} - \mathbf{S}) \mathbf{Y},$$

where $\widetilde{\mathbf{Z}} = (\mathbf{I} - \mathbf{S})\hat{\mathbf{Z}}$, $\mathbf{I}$ is the $n \times n$ identity matrix,

$$\mathbf{S} = \begin{pmatrix} (\mathbf{X}_1^{\mathrm{T}} \quad \mathbf{0}_q^{\mathrm{T}})(\mathbf{D}_{u_1}^{\mathrm{T}} \mathbf{W}_{u_1} \mathbf{D}_{u_1})^{-1} \mathbf{D}_{u_1}^{\mathrm{T}} \mathbf{W}_{u_1} \\ \vdots \\ (\mathbf{X}_n^{\mathrm{T}} \quad \mathbf{0}_q^{\mathrm{T}})(\mathbf{D}_{u_n}^{\mathrm{T}} \mathbf{W}_{u_n} \mathbf{D}_{u_n})^{-1} \mathbf{D}_{u_n}^{\mathrm{T}} \mathbf{W}_{u_n} \end{pmatrix}_{n \times 2q},$$

$$\mathbf{D}_u = \begin{pmatrix} \mathbf{X}_1^{\mathrm{T}} & h^{-1}(U_1 - u)\mathbf{X}_1^{\mathrm{T}} \\ \vdots & \vdots \\ \mathbf{X}_n^{\mathrm{T}} & h^{-1}(U_n - u)\mathbf{X}_n^{\mathrm{T}} \end{pmatrix}_{n \times 2q}$$

and $\mathbf{Y} = (Y_1, \ldots, Y_n)^{\mathrm{T}}$, $\mathbf{W}_u = \mathrm{diag}\{K_h(U_1 - u), \ldots, K_h(U_n - u)\}_{n \times n}$, $\hat{\mathbf{Z}} = (\hat{\mathbf{Z}}_1, \ldots, \hat{\mathbf{Z}}_n)^{\mathrm{T}}$, $\mathbf{0}_q$ is the $q \times 1$ vector with all the entries being zero, $K(\cdot)$ is a kernel function, $h$ is a bandwidth and $K_h(\cdot) = K(\cdot/h)/h$.

We now give a representation of $\hat{\boldsymbol{\Theta}}_n$. This representation can be used to obtain the asymptotic distribution of $\sqrt{n}(\hat{\boldsymbol{\Theta}}_n - \boldsymbol{\Theta})$, which we give in Theorem 2. This result extends the method of Fan and Huang (2005) to a SVCPLM with generated regressors.

Let $\Phi(U) = E(\mathbf{X}\mathbf{Z}^{\mathrm{T}}|U)$, $\Gamma(U) = E(\mathbf{X}\mathbf{X}^{\mathrm{T}}|U)$, $\psi(\mathbf{Z}, \mathbf{X}, U) = \mathbf{Z} - \Phi^{\mathrm{T}}(U) \times \Gamma^{-1}(U)\mathbf{X}$, $B(\mathbf{V}) = E[\{\mathbf{Z} - \Phi^{\mathrm{T}}(U)\Gamma^{-1}(U)\mathbf{X}\}|\mathbf{V}]$ and $\boldsymbol{\Sigma} = E(\mathbf{Z}\mathbf{Z}^{\mathrm{T}}) - E\{\Phi^{\mathrm{T}}(U) \times \Gamma^{-1}(U)\Phi(U)\}$.

THEOREM 1. *Under Assumptions 1–5 in the Appendix, we have*

$$\hat{\mathbf{\Theta}}_n - \mathbf{\Theta} = \mathbf{\Sigma}^{-1}\left[\frac{1}{n}\frac{b^{r+1}}{(r+1)!}\zeta_1^{\mathrm{T}} S_u^{-1} c_p \sum_{i=1}^n \psi(\mathbf{Z}_i, \mathbf{X}_i, U_i)\{\boldsymbol{\xi}^{(r+1)}(\mathbf{V}_i)\}^{\mathrm{T}}\boldsymbol{\beta}_0 \right.$$
$$\left. + \frac{1}{n}\sum_{j=1}^n \Delta(\mathbf{V}_j)\mathbf{e}_j^{\mathrm{T}}\boldsymbol{\beta}_0 + \frac{1}{n}\sum_{i=1}^n \psi(\mathbf{Z}_i, \mathbf{X}_i, U_i)\varepsilon_i\right]$$
$$\times \{1 + o_{\mathrm{P}}(1)\},$$

*where* $\Delta(\mathbf{V}_j) = \frac{1}{n}\sum_{i=1}^n \psi(\mathbf{Z}_i, \mathbf{X}_i, U_i) L_b(\mathbf{V}_j - \mathbf{V}_i)/f_v(\mathbf{V}_i)$.

THEOREM 2. *Let* $nb^{2(r+1)} \to 0$. *Under Assumptions 1–5 in the Appendix,* $\sqrt{n}(\hat{\mathbf{\Theta}}_n - \mathbf{\Theta})$ *converges to a normal distribution with mean zero and covariance matrix* $\mathbf{\Sigma}_1$, *where* $\mathbf{\Sigma}_1 = \mathbf{\Sigma}^{-1}\mathbf{D}\mathbf{\Sigma}^{-1}$, $\mathbf{D} = E[\sigma^2(\mathbf{X}, \mathbf{Z}, U)\{\psi(\mathbf{X}, \mathbf{Z}, U)\}^{\otimes 2}] + E[(\mathbf{e}^{\mathrm{T}}\boldsymbol{\beta})^2\{B(\mathbf{V})\}^{\otimes 2}] + \boldsymbol{\beta}^{\mathrm{T}} E\{E(\mathbf{e}\varepsilon|\mathbf{Z}, \mathbf{X}, U, \mathbf{V})\{B(\mathbf{V})\}^{\otimes 2}\}$.

*Furthermore, if* $\mathbf{e}$ *is independent of* $\varepsilon$ *given* $(\mathbf{Z}, \mathbf{X}, U, \mathbf{V})$, *and* $\varepsilon$ *is independent of* $(\mathbf{Z}, \mathbf{X}, U)$, *the asymptotic covariance can be simplified as* $\mathbf{\Sigma}^{-1}(\sigma^2\mathbf{\Sigma} + E[(\mathbf{e}^{\mathrm{T}}\boldsymbol{\beta})^2\{B(\mathbf{V})\}^{\otimes 2}])\mathbf{\Sigma}^{-1}$. *If* $\mathbf{e}$ *is also independent of* $\mathbf{V}$, *the asymptotic covariance can further be simplified as* $\sigma^2\mathbf{\Sigma}^{-1} + \boldsymbol{\beta}^{\mathrm{T}}\Sigma_e\boldsymbol{\beta}\mathbf{\Sigma}^{-1}E\{B(\mathbf{V})\}^{\otimes 2}\mathbf{\Sigma}^{-1}$.

The proof of Theorem 2 can be completed by using Theorem 1. We omit the details.

The asymptotic variance has a similar structure to that of Das (2005). The first term of asymptotic variance can be viewed as the variance from the first stage estimation without measurement error/missing data, the second one is the variance of the second stage for estimating unobserved variables and the third one is the covariance of two-stage estimators. If $\mathbf{e} = 0$ in (1.2), that is, the covariate can be exactly observed, the variance of $\hat{\mathbf{\Theta}}_n$ is the same as that of Fan and Huang (2005). To achieve the root-$n$ estimator of $\mathbf{\Theta}$, Theorem 2 indicates that undersmoothing is required in estimating $\boldsymbol{\xi}(v)$ and the optimal bandwidth does not satisfy the condition of Theorem 2.

EXAMPLE 1 (cont.). Let $\hat{\boldsymbol{\beta}}_n$ be the estimator of $\boldsymbol{\beta}$ in (1.3). Assume $n_0/n \to \lambda$. Checking the conditions of Theorem 2, we can conclude that $\sqrt{n}(\hat{\boldsymbol{\beta}}_n - \boldsymbol{\beta}_0) \xrightarrow{\mathcal{L}} N(0, \mathbf{\Sigma}_\star)$, where $\mathbf{\Sigma}_\star = \mathbf{\Sigma}^{-1}(\sigma^2 + \lambda\boldsymbol{\beta}^{\mathrm{T}} E[E\{\mathbf{Z} - E(\mathbf{Z}|U)|\mathbf{V}\}]^{\otimes 2}\boldsymbol{\beta})$ and $\mathbf{\Sigma} = E[\{\boldsymbol{\xi} - E(\boldsymbol{\xi}|U)\}^{\otimes 2}]$.

EXAMPLE 2 (cont.). For the de-noised models introduced in Section 1, we apply Theorem 2 to derive the asymptotic distribution of the estimator, $\hat{\mathbf{\Theta}} = (\hat{\boldsymbol{\beta}}^{\mathrm{T}}, \hat{\theta}^{\mathrm{T}})^{\mathrm{T}}$, given by Cui, He and Zhu (2002), and obtain that $\sqrt{n}(\hat{\mathbf{\Theta}} - \mathbf{\Theta}) \xrightarrow{\mathcal{L}} N\{0, \mathbf{\Sigma}^{-1}(\sigma^2 + \boldsymbol{\beta}^{\mathrm{T}}\Sigma_e\boldsymbol{\beta})\}$.

The asymptotic covariance of $\hat{\boldsymbol{\Theta}}_n$ can be consistently estimated by $\hat{\boldsymbol{\Sigma}}_n = n\hat{\boldsymbol{\Sigma}}^{-1}\hat{\sigma}^2 + \hat{\boldsymbol{\Sigma}}^{-1}\hat{\mathbf{Q}}\hat{\boldsymbol{\Sigma}}^{-1}$, where $\hat{\boldsymbol{\Sigma}}^{-1} = \{(\tilde{\mathbf{Z}}^T\tilde{\mathbf{Z}})^{-1}\tilde{\mathbf{Z}}^T(\mathbf{I}-\mathbf{S})^T\}^{\otimes 2}$, $\hat{\mathbf{Q}} = \frac{1}{n}\sum_{i=1}^n(\boldsymbol{\eta}_i - \hat{\mathbf{Z}}_i)^T\hat{\boldsymbol{\Theta}}_n\{\hat{\mathbf{B}}(\mathbf{V}_i)\}^{\otimes 2}$, $\hat{\sigma}^2 = \frac{1}{n}\sum_{i=1}^n\{Y_i - \hat{\alpha}(U_i)\mathbf{X}_i - \hat{\boldsymbol{\Theta}}^T\hat{\mathbf{Z}}\}^2$, $\hat{\mathbf{B}}(\mathbf{v}) = \hat{\mathbf{Z}} - \hat{E}\{\Phi^T(U)\Gamma^{-1}(U)\mathbf{X}|\mathbf{V}=\mathbf{v}\}$ and $\hat{E}\{\Phi^T(U)\Gamma^{-1}(U)\mathbf{X}|\mathbf{V}=\mathbf{v}\}$ is a nonparametric regression estimator of $\Phi^T(U)\Gamma^{-1}(U)\mathbf{X}$ on $\mathbf{V}$. $\hat{\alpha}(\cdot)$ will be given in the next section.

Generally $\hat{\boldsymbol{\Sigma}}_n$ is difficult to calculate. However, implementation will become simpler in some cases. For example, in the errors-in-variables model with validation data, a direct simplification yields $B(\mathbf{V}) = \mathbf{Z} - \Phi^T(U)\Gamma^{-1}(U)\mathbf{X}$, $\mathbf{D} = \{\boldsymbol{\beta}^T E(\mathbf{ee}^T|\mathbf{V})\boldsymbol{\beta}\}\boldsymbol{\Sigma}$ and the asymptotic covariance matrix equals $\boldsymbol{\Sigma}^{-1}\{\sigma^2 + \lambda\boldsymbol{\beta}^T E(\mathbf{ee}^T|\mathbf{V})\boldsymbol{\beta}\}$. This matrix can be estimated by a standard sandwich procedure. The similar situation also applies for the asymptotic covariance matrix, $\boldsymbol{\Sigma}^{-1}\{\sigma^2 + \boldsymbol{\beta}^T E(\mathbf{ee}^T|\mathbf{V})\boldsymbol{\beta}\}$, of the de-noise model.

2.3. *Estimation of the nonparametric components.* After obtaining estimates $\hat{\boldsymbol{\Theta}}_n$, we can estimate $a_j(u)$ and $b_j(u)$ for $j = 1, \ldots, q$, and then $\boldsymbol{\alpha}_j(u)$. Write $\boldsymbol{\Psi}(u) = \{a_1(u), \ldots, a_q(u), b_1(u), \ldots, b_q(u)\}^T$. An estimator of the nonparametric components $\boldsymbol{\Psi}(u)$ is defined as

$$(2.4) \qquad \hat{\boldsymbol{\Psi}}(u) = \mathbf{H}^{-1}(\mathbf{D}_u^T\mathbf{W}_u\mathbf{D}_u)^{-1}\mathbf{D}_u^T\mathbf{W}_u(\mathbf{Y} - \hat{\mathbf{Z}}\hat{\boldsymbol{\Theta}}_n).$$

Correspondingly, $\mathbf{a}(u)$ is estimated by $\hat{\mathbf{a}}(u) = (\mathbf{I}_q, \mathbf{0}_q)(\mathbf{D}_u^T\mathbf{W}_u\mathbf{D}_u)^{-1}\mathbf{D}_u^T\mathbf{W}_u(\mathbf{Y} - \hat{\mathbf{Z}}\hat{\boldsymbol{\Theta}}_n)$, where $\mathbf{I}_q$ is the $q \times q$ identity matrix, $\mathbf{H} = \text{diag}(1, h) \otimes \mathbf{I}_q$. We have the following asymptotic representation for the resulting estimator:

THEOREM 3. *Under Assumptions 1–5 given in the [Appendix](), we have*
$$\sqrt{nh}\mathbf{H}\{\hat{\boldsymbol{\Psi}}(u_0) - \boldsymbol{\Psi}(u_0)\}$$
$$= \frac{n^{1/2}h^{5/2}}{2(\mu_2 - \mu_1^2)}\begin{pmatrix}(\mu_2^2 - \mu_1\mu_3)\\(\mu_3 - \mu_1\mu_2)\end{pmatrix}\boldsymbol{\alpha}''(u_0)$$
$$- \frac{\sqrt{nh}b^{r+1}}{(r+1)!}\zeta_1^T S_u^{-1} c_p \begin{pmatrix}\Gamma^{-1}(u_0)E[\mathbf{X}\{\boldsymbol{\xi}^{(r+1)}(\mathbf{V})\}^T\boldsymbol{\beta}_0|U = u_0]\\0\end{pmatrix}$$
$$+ o(n^{1/2}h^{5/2} + n^{1/2}h^{1/2}b^{r+1})$$
$$+ \frac{\sqrt{nh}\Gamma^{-1}(u)}{nf_u(u)(\mu_2 - \mu_1^2)}\sum_{i=1}^n K_h(U_i - u)\{\mathbf{X}_i\varepsilon_i + E(\mathbf{X}_i|\mathbf{V}_i)\mathbf{e}_i^T\boldsymbol{\beta}\}$$
$$\otimes \begin{pmatrix}\mu_2 - \mu_1(U_i - u)/h\\(U_i - u)/h - \mu_1\end{pmatrix}\{1 + o_P(1)\}.$$

Based on this representation, we can derive the asymptotic normality of the proposed nonparametric estimators of the varying coefficient functions. The proof is straightforward but tedious. We omit the details.

For notational simplicity, we assume that $\varepsilon$ is independent of $(\mathbf{Z}, \mathbf{X}, U)$ and $\mathbf{e}$ is independent of $(\mathbf{V}, U)$ in the remaining part of this paper.

THEOREM 4. *Under Assumptions 1–5, we have*

$$\sqrt{nh}\bigg[\mathbf{H}\{\hat{\boldsymbol{\Psi}}(u_0) - \boldsymbol{\Psi}(u_0)\} - \frac{h^2}{2(\mu_2 - \mu_1^2)}\begin{pmatrix}(\mu_2^2 - \mu_1\mu_3)\boldsymbol{\alpha}''(u_0)\\(\mu_3 - \mu_1\mu_2)\boldsymbol{\alpha}''(u_0)\end{pmatrix}$$

$$- \frac{b^{r+1}}{(r+1)!}\zeta_1^{\mathrm{T}} S_u^{-1} c_p$$

$$\times \begin{pmatrix}\Gamma^{-1}(u_0)E[\mathbf{X}\{\boldsymbol{\xi}^{(r+1)}(\mathbf{V})\}^{\mathrm{T}}\boldsymbol{\beta}_0|U = u_0]\\0\end{pmatrix} + o(h^2 + b^{r+1})\bigg]$$

$$\xrightarrow{\mathcal{L}} N(0, \boldsymbol{\Sigma}_2),$$

*as* $n \to \infty$, *where* $\boldsymbol{\Sigma}_2 = f_u^{-1}(u_0)\{\sigma^2\Gamma^{-1}(u_0) + \Gamma^{-1}(u_0)\boldsymbol{\Sigma}_1^*\Gamma^{-1}(u_0)\} \otimes \mathcal{G}$,

$$\mathcal{G} = \frac{1}{(\mu_2 - \mu_1^2)^2}$$
$$\times \begin{pmatrix}\mu_2^2\nu_0 - 2\mu_1\mu_2\nu_1 + \mu_1^2\nu_2 & (\mu_1^2 + \mu_2)\nu_1 - \mu_1\mu_2\nu_0 - \mu_1\nu_2\\(\mu_1^2 + \mu_2)\nu_1 - \mu_1\mu_2\nu_0 - \mu_1\nu_2 & \nu_2 - \mu_1(2\nu_1 + \mu_1\nu_0)\end{pmatrix},$$

$\boldsymbol{\Sigma}_1^* = \boldsymbol{\beta}^{\mathrm{T}}\Sigma_e\boldsymbol{\beta}\Lambda(u_0)$, $\Lambda(u_0) = (E[\{E(\mathbf{X}|\mathbf{V})\}|U = u_0])^{\otimes 2}$, $q_0 = \mu_2/(\mu_2 - \mu_1)$, $q_1 = -\mu_1/(\mu_2 - \mu_1^2)$.

*Furthermore, if* $nhb^{2r+2} \to 0$, *then*

$$\sqrt{nh}\bigg\{\hat{\boldsymbol{\alpha}}(u) - \boldsymbol{\alpha}(u) - \frac{h^2}{2}\frac{\mu_2^2 - \mu_1\mu_3}{\mu_2 - \mu_1^2}\boldsymbol{\alpha}''(u) + o(h^2 + b^{r+1})\bigg\} \xrightarrow{L} N(0, \boldsymbol{\Sigma}_2^*),$$

*where* $\boldsymbol{\Sigma}_2^* = \sigma^2(q_0^2\nu_0 + 2q_0q_1\nu_1 + q_1^2\nu_2)\{\Gamma^{-1}(u_0) + \Gamma^{-1}(u_0)\boldsymbol{\Sigma}_1^*\Gamma^{-1}(u_0)\}/f_u(u)$.

The first term of $\boldsymbol{\Sigma}_2$ is the asymptotic covariance of the usual profile likelihood estimator of Cai, Fan and Li (2000), when $\xi_j$ is observed. The second term is attributed to calibrating the error-prone covariates. In the error-in-variable model with validation data, if $\mathbf{X}$ is independent of $\mathbf{V}$ and $E(\mathbf{X}) = 0$, the measurement errors have no impact on the effect of the covariance $\boldsymbol{\Sigma}_2$. Theorem 4 also indicates that if $n^{1/2}\max(h^{5/2}, b^{r+1}) \to 0$, the bias of $\hat{\alpha}(u)$ tends to zero and $\hat{\alpha}(u)$ is asymptotically normally distributed with rate $(nh)^{1/2}$.

After obtaining $\hat{\boldsymbol{\Theta}}_n$ and $\hat{\alpha}(u)$, one can easily give an estimator of the variance $\sigma^2$ of the error $\varepsilon$:

$$\hat{\sigma}_n^2 = \frac{1}{n}\sum_{i=1}^n \{Y_i - \hat{\boldsymbol{\beta}}^{\mathrm{T}}\hat{\boldsymbol{\xi}}_n^{\mathrm{T}}(V_i) - \hat{\boldsymbol{\theta}}_n^{\mathrm{T}}\mathbf{W}_i - \hat{\alpha}^{\mathrm{T}}(U_i)\mathbf{X}_i\}^2.$$

In our simulation, a simple version of $\hat{\sigma}_n^2$ is used. Note that $\mathbf{S}$ depends only on the observations $\{(U_i, \mathbf{X}_i)\}_{i=1}^n$, and we can derive a "synthetic linear model," that is, $\mathbf{Y} - \mathbf{Z\Theta} = \mathbf{M} + \boldsymbol{\varepsilon}$, where $\mathbf{M} = \boldsymbol{\alpha}^{\mathrm{T}}(U)\mathbf{X}$. A straightforward derivation yields $(\mathbf{I} - \mathbf{S})\mathbf{Y} = (\mathbf{I} - \mathbf{S})\mathbf{Z\Theta} + (\mathbf{I} - \mathbf{S})\boldsymbol{\varepsilon}$. Standard regression gives the least-square estimates $\hat{\mathbf{\Theta}}$ and then $\hat{\mathbf{M}} = \mathbf{S}(\mathbf{Y} - \mathbf{Z}\hat{\mathbf{\Theta}})$. Note that $\mathbf{Z}$ is not always observed. Replacing $\mathbf{Z}$ by its estimates, we obtain a consistent estimator $\hat{\mathbf{M}}$ of $\mathbf{M}$; that is, $\hat{\mathbf{M}} = \mathbf{S}(\mathbf{Y} - \hat{\mathbf{Z}}\hat{\mathbf{\Theta}})$. A consistent estimator $\sigma^2$ may be defined as $\hat{\sigma}_n^2 = \frac{1}{n}\sum_{i=1}^n (Y_i - \hat{\mathbf{\Theta}}^{\mathrm{T}}\hat{\mathbf{Z}}_i - \hat{M}_i)^2$, where $\hat{M}_i$ is the $i$th element of $\hat{\mathbf{M}}$.

2.4. *Bandwidth selection.* The proposed procedure involves two bandwidths, $h$ and $b$, to be selected. To derive asymptotic distributions of the proposed estimators, we theoretically impose the rates of convergence for the bandwidths. It is worthwhile to point out that undersmoothing is necessary when we estimate $\xi$ and the optimal bandwidth for $b$ is then violated.

As mentioned before, the optimal bandwidth for $b$ cannot be obtained because undersmoothing the nonparametric estimators of the covariates is necessary. The consequence of undersmoothing $\boldsymbol{\xi}$ is that the bias is kept small and precludes the optimal bandwidth for $b$. The asymptotic variances of the proposed estimators for constant coefficients depend on neither the bandwidth nor the kernel function. Hence, we can use the similar method of mixture of higher-order theoretical expansions, proposed by Sepanski, Knicherbocker and Carroll (1994) or the typical curves approach by Brookmeyer and Liao (1992) to select the bandwidth $b$. As done by Sepanski, Knickerbocker and Carroll (1994), the suitable bandwidth is $b = Cn^{-1/3}$, where $C$ is a constant depending on unknown function $\boldsymbol{\xi}(v)$ and its twice derivatives. In practice, one can use a plug-in rule to estimate the constant $C$. A useful and simple candidate $C$ is $\hat{\sigma}_V$, the sample deviation of $V$. This method is fairly effective and easy to implement. In our simulation example, the bandwidth is $b = \hat{\sigma}_v n^{-1/3}$. Based on the asymptotic analysis and empirical experience for the fixed time case (i.e., de-noise models), we suggest a simple rule of thumb as follows: The smoothing parameter $b$ is so chosen that intervals of size $2b$ would contain around 5 points for $n$ up to 100 and between $8^{-1}n^{1/3}$ and $4^{-1}n^{1/3}$ points for larger $n$.

We use the "leave one sample out" method to select the bandwidth $h$. This method has been widely applied in practice; for example, Cai, Fan and Li (2000) and Fan and Huang (2005). We define the cross-validation score for $h$ as $CV(h) = n^{-1}\sum_{i=1}^n \{Y_i - \hat{\alpha}_{h,-i}^{\mathrm{T}}(U_i)\mathbf{X}_i - \hat{\mathbf{\Theta}}_{n,-i}^{\mathrm{T}}\hat{\mathbf{Z}}_i\}^2$, where $\hat{\mathbf{\Theta}}_{n,-i}$ is the estimated profile least-square-based estimator defined by (2.3), computed from the data with measurements of the $i$th observation deleted, and $\hat{\boldsymbol{\alpha}}_{h,-i}(\cdot)$ is the estimator defined in (2.4) with $\hat{\mathbf{\Theta}}_n$ replaced by $\hat{\mathbf{\Theta}}_{n,-i}$. The likelihood cross-validation smoothing parameter $h_{cv}$ is the minimizer of $CV(h)$. That is, $h_{cv} = \arg\min_h CV(h)$.

### 3. Tests for parametric and nonparametric components.

3.1. *Test for parametric components.* An interesting question is to consider the following hypothesis:

(3.1) $$H_0: \mathbf{A}\boldsymbol{\Theta} = 0 \quad \text{VS} \quad H_1: \mathbf{A}\boldsymbol{\Theta} \neq 0,$$

where $\mathbf{A}$ is a given $l \times p$ full rank matrix.

Let $\hat{\boldsymbol{\Theta}}_0 = (\hat{\boldsymbol{\beta}}_0^{\mathrm{T}}, \hat{\boldsymbol{\theta}}_0^{\mathrm{T}})^{\mathrm{T}}$ be the estimators of $\boldsymbol{\Theta}$ and $\hat{\boldsymbol{\alpha}}_0(\cdot)$ be the estimator of $\alpha(u)$ under the null hypothesis. Denote $RSS_0 = \sum_{i=1}^n \{Y_i - \hat{\boldsymbol{\beta}}_0^{\mathrm{T}} \hat{\boldsymbol{\xi}}_i - \hat{\boldsymbol{\theta}}_0^{\mathrm{T}} \mathbf{W}_i - \hat{\boldsymbol{\alpha}}_0^{\mathrm{T}}(U_i)\mathbf{X}_i\}^2$. $RSS_0$ can be further expressed as $\sum_{i=1}^n \{Y_i - \hat{\boldsymbol{\beta}}_0^{\mathrm{T}} \hat{\boldsymbol{\xi}}_i - \hat{\boldsymbol{\theta}}_0^{\mathrm{T}} \mathbf{W}_i - \mathbf{S}(Y - \hat{\mathbf{Z}}\hat{\boldsymbol{\Theta}}_0)\}^2$, where $\hat{\boldsymbol{\Theta}}_0 = \hat{\boldsymbol{\Theta}} - (\tilde{\mathbf{Z}}^{\mathrm{T}}\tilde{\mathbf{Z}})^{-1}\mathbf{A}^{\mathrm{T}}\{\mathbf{A}(\tilde{\mathbf{Z}}^{\mathrm{T}}\tilde{\mathbf{Z}})^{-1}\mathbf{A}^{\mathrm{T}}\}^{-1}\mathbf{A}\hat{\boldsymbol{\Theta}}$, and $\hat{\boldsymbol{\Theta}} = (\tilde{\mathbf{Z}}^{\mathrm{T}}\tilde{\mathbf{Z}})^{-1}\tilde{\mathbf{Z}}^{\mathrm{T}}\tilde{\mathbf{Y}}$, an estimator of $\boldsymbol{\Theta}$ without the restriction, with $\tilde{\mathbf{Z}} = (\mathbf{I} - \mathbf{S})\hat{\mathbf{Z}}$ and $\tilde{\mathbf{Y}} = (\mathbf{I} - \mathbf{S})\hat{\mathbf{Y}}$.

Similar, let $\hat{\boldsymbol{\Theta}}_1 = (\hat{\boldsymbol{\beta}}_1^{\mathrm{T}}, \hat{\boldsymbol{\theta}}_1^{\mathrm{T}})^{\mathrm{T}}$ and $\hat{\boldsymbol{\alpha}}_1(\cdot)$ be the estimators of $\boldsymbol{\Theta}$ and $\alpha(\cdot)$ under the alternative hypothesis, respectively. Denote $RSS_1 = \sum_{i=1}^n \{Y_i - \hat{\boldsymbol{\beta}}_1^{\mathrm{T}} \hat{\boldsymbol{\xi}}_i - \hat{\boldsymbol{\theta}}_1^{\mathrm{T}} \mathbf{W}_i - \hat{\boldsymbol{\alpha}}_1^{\mathrm{T}}(U_i)\mathbf{X}_i\}^2$, which can be expressed as $\sum_{i=1}^n \{Y_i - \hat{\boldsymbol{\beta}}_1^{\mathrm{T}} \hat{\boldsymbol{\xi}}_i - \hat{\boldsymbol{\theta}}_1^{\mathrm{T}} \mathbf{W}_i - S(Y - \hat{\mathbf{Z}}\hat{\boldsymbol{\Theta}}_1)\}^2$. Following Fan and Huang (2005), we define a profile least-square-based ratio test by

$$T_n = \frac{n}{2}(RSS_0 - RSS_1)/RSS_1.$$

Under their set-up, Fan and Huang (2005) showed that statistic $T_n$ is the profile likelihood ratio when the error distribution is normally distributed. In the present situation, because of the effect of measurement error on variables, no central $\mathcal{X}^2$-distribution similar to that of Fan and Huang (2005) is available. However, we can still prove that $2T_n$ has the asymptotic noncentral $\chi^2$ distribution under the alternative hypothesis of (3.1), which we summarize in the following theorem.

THEOREM 5. *Suppose that Assumptions 1–5 in the Appendix are satisfied and $nb^{2r+2} \to 0$, as $n \to \infty$. Under the alternative hypothesis of (3.1),*

$$2T_n - n\sigma^{-2}\boldsymbol{\Theta}^{\mathrm{T}}\mathbf{A}^{\mathrm{T}}(\mathbf{A}\boldsymbol{\Sigma}^{-1}\mathbf{A}^{\mathrm{T}})^{-1}\mathbf{A}\boldsymbol{\Theta} \xrightarrow{L} \sum_{i=1}^l \omega_i \chi_{i1}^2$$

*where $\omega_i$ for $1 \leq i \leq l$ are the eigenvalues of $(\sigma^2 \mathbf{A}\boldsymbol{\Sigma}^{-1}\mathbf{A}^{\mathrm{T}})^{-1}(\mathbf{A}\boldsymbol{\Sigma}_1^{-1}\mathbf{A}^{\mathrm{T}})$ and $\chi_{i1}^2$ is the central $\chi^2$ distribution with 1 degree of freedom. Furthermore, let $\hat{\boldsymbol{\Sigma}}_1$ and $\hat{\boldsymbol{\Sigma}}$ be the consistent estimators of $\boldsymbol{\Sigma}_1$ and $\boldsymbol{\Sigma}$, respectively. Then $2\varrho_n T_n \xrightarrow{L} \chi_{(l)}^2(\lambda)$, where $\varrho_n = l/\operatorname{tr}\{(\sigma^2 \mathbf{A}\hat{\boldsymbol{\Sigma}}^{-1}\mathbf{A}^{\mathrm{T}})^{-1}(\mathbf{A}\hat{\boldsymbol{\Sigma}}_1^{-1}\mathbf{A}^{\mathrm{T}})\}$, $\chi_{(l)}^2(\lambda)$ is the noncentral $\chi^2$ distribution with $l$ degree of freedom, and the noncentral parameter $\lambda = \sigma^{-2}\varrho \lim_{n\to\infty} n\boldsymbol{\Theta}^{\mathrm{T}}\mathbf{A}^{\mathrm{T}}(\mathbf{A}\boldsymbol{\Sigma}^{-1}\mathbf{A}^{\mathrm{T}})^{-1}\mathbf{A}\boldsymbol{\Theta}$ with $\varrho = l/\operatorname{tr}\{(\sigma^2 \mathbf{A}\boldsymbol{\Sigma}^{-1}\mathbf{A}^{\mathrm{T}})^{-1} \times (\mathbf{A}\boldsymbol{\Sigma}_1^{-1}\mathbf{A}^{\mathrm{T}})\}$.*

In a similar way, we may construct the Wald test for hypothesis (3.1) as $W_n = \hat{\boldsymbol{\Theta}}^{\mathrm{T}} \mathbf{A}^{\mathrm{T}} (\mathbf{A} \hat{\boldsymbol{\Sigma}}_1 \mathbf{A}^{\mathrm{T}})^{-1} \mathbf{A} \hat{\boldsymbol{\Theta}}$, and demonstrate that $W_n$ and $2\varrho_n T_n$ have the same asymptotic distribution under the alternative hypothesis. These properties can therefore be used to calculate the power of the proposed tests.

EXAMPLE 3 (cont.).  Generalize (1.6) to a more flexible model:

$$Y_t = \boldsymbol{\beta}^{\mathrm{T}} E(\boldsymbol{\eta}|\mathbf{V}_t) + \boldsymbol{\zeta}^{\mathrm{T}} \boldsymbol{\eta} + \boldsymbol{\gamma}^{\mathrm{T}} \mathbf{S}_t + \alpha(U_t)\mathbf{X}_t + \varepsilon_t.$$

Write $\boldsymbol{\Theta} = (\boldsymbol{\beta}^{\mathrm{T}}, \boldsymbol{\zeta}^{\mathrm{T}}, \boldsymbol{\gamma}^{\mathrm{T}})^{\mathrm{T}}$ and $\mathbf{Z} = \{E(\boldsymbol{\eta}^{\mathrm{T}}|\mathbf{V}), \boldsymbol{\eta}^{\mathrm{T}}, \mathbf{S}_t^{\mathrm{T}}\}^{\mathrm{T}}$. The hypothesis (1.7) is equivalent to

$$(3.2) \qquad \mathbf{A}\boldsymbol{\Theta} = 0 \quad \text{VS} \quad H_1 : \mathbf{A}\boldsymbol{\Theta} \neq 0,$$

where $\mathbf{A} = (\mathbf{1}_{p_1}, -\mathbf{1}_{p_1}, 0)$, $\mathbf{1}_{p_1}$ is $p_1$-variate vector with all entries 1. This is an expression of (3.1). As a consequence, the proposed profile least-square-based ratio test and Wald test can be applied to test this hypothesis.

For hypothesis (3.2), one may also propose a Wald-type statistic: $W_n(h) = \hat{\boldsymbol{\Theta}}^{\mathrm{T}} \mathbf{A}^{\mathrm{T}} (\mathbf{A} \hat{\boldsymbol{\Sigma}}_h \mathbf{A}^{\mathrm{T}})^{-1} \mathbf{A} \hat{\boldsymbol{\Theta}}$, where $\hat{\boldsymbol{\Sigma}}_h = \hat{\Sigma}^{-1}(\hat{\sigma}^2 + \hat{\boldsymbol{\beta}}^{\mathrm{T}} \hat{\Sigma}_e \hat{\boldsymbol{\beta}})$. It can be proved that $2\varrho_n T_n$ and $W_n$ have the same asymptotic $\mathcal{X}^2$ distribution.

3.2. *Tests for the nonparametric part and wild bootstrap version.* It is also of interest to check whether the varying-coefficient functions $\alpha(u)$ in (1.2) are parametric functions. Specifically speaking, we consider the following hypothesis:

$$H_0 : \alpha_i(U) = \alpha_i(U, \boldsymbol{\gamma}) \quad \text{VS} \quad H_1 : \alpha_i(U) \neq \alpha_i(U, \boldsymbol{\gamma}), \qquad i = 1, 2, \ldots, q,$$

where $\boldsymbol{\gamma}$ is an unknown vector, $\alpha_i(\cdot, \cdot)$ is a known function and $i = 1, 2, \ldots, q$.

For simplicity of presentation, we test the homogeneity:

$$H_0 : \alpha_1(U) = \alpha_1, \ldots, \alpha_q(U) = \alpha_q.$$

Let $\tilde{\alpha}_1, \ldots, \tilde{\alpha}_q$ and $\tilde{\boldsymbol{\Theta}}$ be the profile estimator under $H_0$. The weighted residual sum of squares under $H_0$ is $RSS(H_0) = \sum_{i=1}^n w_i (Y_i - \sum_{j=1}^q \tilde{\alpha}_j X_{ij} - \tilde{\boldsymbol{\Theta}}^{\mathrm{T}} \hat{\mathbf{Z}}_i)^2$, where $w_i(\cdot)$ are weighted functions such that $\sum_{i=1}^n w_i = 1$, and $w_i \geq 0$. In general, the weight function $w$ has a compact support, designed to reduce the boundary effects on the test statistics. When $\sigma^2(\mathbf{Z}, \mathbf{X}, U) = v(\mathbf{Z}, \mathbf{X}, U)\sigma^2$ for some known function $v(\mathbf{Z}, \mathbf{X}, U)$, we may choose $w_i = v^{-1}(\mathbf{Z}_i, \mathbf{X}_i, U_i)$. See Fan, Zhang and Zhang (2001) and Fan and Jiang (2007) for a similar argument.

Under the general alternative that all the varying-coefficient functions are allowed to be varying of random variable $U$, we use the local likelihood method to obtain estimator $\hat{\boldsymbol{\beta}}$ and $\hat{\boldsymbol{\alpha}}(U)$. Therefore, the corresponding weighted residual sum of squares is

$$RSS(H_1) = \sum_{i=1}^{n} w_i \left\{ Y_i - \sum_{j=1}^{q} \hat{\alpha}_j(U_i) X_{ij} - \hat{\boldsymbol{\Theta}}^{\mathrm{T}} \hat{\mathbf{Z}}_i \right\}^2.$$

In a similar way to that used in Section 3.1, we propose a generalized likelihood ratio (GLR) statistic: $T_{\mathrm{GLR}} = \{RSS(H_0) - RSS(H_1)\}/RSS(H_1)$. Under mild assumptions, one can derive the asymptotic distribution of $T_{\mathrm{GLR}}$. This distribution can be used to gain the empirical level. See Fan, Zhang and Zhang (2001) for a related discussion.

These arguments can be applied to the following partially parametric null hypothesis: $H_0: \alpha_1(U) = \alpha_1, \ldots, \alpha_l(U) = \alpha_r, r < q$. The difference is only the definition of $RSS(H_0)$, for which we use the profile likelihood procedure to estimate the constant coefficient $\alpha_i$, $i = 1, 2, \ldots, r$ and $\boldsymbol{\Theta}$, and use the profile linear procedure to estimate the nonparametric component $\alpha_i(\cdot)$, $i = r+1, \ldots, q$ under the null hypothesis.

Although the asymptotic level of $T_{\mathrm{GLR}}$ is available, $T_{\mathrm{GLR}}$ may not perform well when sample sizes are small. For this reason and for practical purposes, we suggest using a bootstrap procedure. To be specific, let $\hat{\varepsilon}_i = Y_i - \hat{\boldsymbol{\Theta}}^{\mathrm{T}} \hat{\mathbf{Z}}_i - \hat{\boldsymbol{\alpha}}^{\mathrm{T}}(U_i) \mathbf{X}_i$ be the residuals based on estimators (2.3) and (2.4) for parametric and nonparametric parts, respectively. We use the Wild bootstrap [Wu (1986), Härdle and Mammen (1993)] method to calculate the critical values for test $T_{\mathrm{GLR}}$. Let $\tau$ be a random variable with a distribution function $F(\cdot)$ such that $E\tau = 0$, $E\tau^2 = 1$ and $E|\tau|^3 < \infty$. We generate the bootstrap residual $\varepsilon_i^* = \hat{\varepsilon}_i \tau_i$, where $\tau_i$ is independent of $\hat{\varepsilon}_i$. Define bootstrap version $T_{\mathrm{GLR}}^*$ like $T_{\mathrm{GLR}}$ based on the bootstrap sample $(Y_i^*, \mathbf{X}_i, \hat{\mathbf{Z}}_i, U_i)$, where $Y_i^* = \hat{\boldsymbol{\Theta}} \hat{\mathbf{Z}}_i + \hat{\boldsymbol{\alpha}}(U_i) \mathbf{X}_i + \varepsilon_i^*$ for $i = 1, 2, \ldots, n$. On a basis of the distribution of $T_{\mathrm{GLR}}^*$, we have the $(1-\alpha)$ quantile $t_{1-\alpha}^*$ and reject the parametric hypothesis if $T_{\mathrm{GLR}} > t_{1-\alpha}^*$.

## 4. Numerical examples.

4.1. *Performance of the proposed estimators.* In this section, we conducted simulation experiments to illustrate the finite sample performances of the proposed estimators and tests. Our simulated data were generated from the following model:

(4.1) $\begin{cases} Y = \beta_1 \xi + \beta_2 W_1 + \beta_3 W_2 + \alpha_1(U) X_1 + \alpha_2(U) X_2 + \varepsilon, \\ \xi = \xi(V), \ \eta = \xi(V) + e. \end{cases}$

$W_1$ and $W_2$ are bivariate normal with marginal mean zero, marginal variance 1 and correlation $1/\sqrt{5}$, while $X_1$ and $X_2$ are independent and normal with mean zero and variance 0.8. The unobserved covariate $\xi$ is related to auxiliary variable $(\eta, V)$ through $\xi(V) = 3V - 2\cos(4\pi V)$ and $\eta = \xi(V) + e$. $V$ is a uniform random variable on $[0,1]$ and $U$ is a uniform random variable on $[0,3]$. The errors $\varepsilon$ and $e$ are independent of each other and normal variables with mean 0 and variances $\sigma_\varepsilon^2$ and $\sigma_e^2$, respectively. The varying-coefficient functions are

$$\alpha_1(U) = \exp(-U^2) + \sin(\pi U) \quad \text{or} \tag{4.2}$$

$$\widetilde{\alpha}_1(U, \varrho) = m + \varrho\{\alpha_1(U) - m\}, \tag{4.3}$$

$$\alpha_2(U) = \tfrac{1}{2}U^2 - \cos(2\pi U), \tag{4.4}$$

where $m = \int_0^3 \alpha_1(t)\,dt/3$, and $\varrho$ is chosen one from the set $\{0.0, 0.2, 0.5, 0.7, 1.0\}$.

The sample size was 100. We generated 500 data sets in each case, applying to each simulated sample the bootstrap test proposed for the parametric part based on 500 bootstrap repetitions. The Gaussian kernel has been used in this example. The optimal bandwidth $h$ was chosen by the leave one out cross-validation method described in Section 2.4 and the bandwidth $b$ was selected as $b = \sigma_v n^{-1/3}$, where $\sigma_v$ is the sample deviation of $V$.

We consider four scenarios. In the first three scenarios $\sigma_\varepsilon^2 = 1$ and $\sigma_e^2 = 2$.

(i) $\beta = (0, c-1, 1)^{\mathrm{T}}$ for $c \in \{0, 0.1, 0.2, 0.25, 0.5, 0.7, 1.0\}$ and $\alpha_1(u)$ and $\alpha_2(u)$ are given in (4.2) and (4.4);

(ii) $\beta = (0, -0.8, 1)^{\mathrm{T}}$ and $\alpha_1(u)$ and $\alpha_2(u)$ are given in (4.3) and (4.4) with $\varrho \in \{0.0, 0.2, 0.5, 0.7, 1.0\}$;

(iii) $\beta = (0.2, -1, 1)^{\mathrm{T}}$ and $\alpha_1(u)$ and $\alpha_2(u)$ are the same as in (ii);

(iv) The setting is the same as that of (iii). But the signal-noise ratio $(r = \sigma_\xi^2/(\sigma_\xi^2 + \sigma_e^2))$ varies from 0.3 to 0.8 by 0.1.

The corresponding results are presented in Tables 1–4, in which we display the estimated values and associated standard errors, standard derivations, and coverage probabilities based on the benchmark estimator (i.e., all covariates measured exactly), the proposed estimator and the naive estimator ($\eta_i$ directly used as the covariates). We summarize our findings as follows:

When $\beta_1 = 0$ [scenario (i) and (ii)], all estimates are close to the true values regardless of the nonparametric functions $\alpha_1(u)$ and $\alpha_2(u)$. The differences among the estimated values based on three methods are slight and can be ignored. However, when $\beta_1 = 0.2$, the estimates of $\beta_1$ based on the naive method have severe biases and the associated coverage probabilities are also substantially smaller than 0.95. These biases were not improved when the sample size was increased (not listed here). But the proposed estimation procedure performs well. On the other hand, the estimates of $\beta_2$ and $\beta_3$ are similar based on the three methods. From Table 4, we can see that the naive estimator of $\beta_1$ has zero coverage probabilities when $r = 0.3$, while the proposed estimator has fairly reasonable coverage probabilities. With an increase of $r$, it is readily seen that coverage probabilities of the proposed estimator are closer to the nominal level, which indicates the proposed method is promising.

TABLE 1
*Results of simulation study for scenario* (i)

| $\varrho$ | | $\beta_1$ | | | | $\beta_2$ | | | | $\beta_3$ | | | |
|---|---|---|---|---|---|---|---|---|---|---|---|---|---|
| | | **Est.** | **SE** | **SD** | **COV** | **Est.** | **SE** | **SD** | **COV** | **Est.** | **SE** | **SD** | **COV** |
| 0 | B | −0.000 | 0.030 | 0.027 | 0.912 | −0.990 | 0.133 | 0.125 | 0.938 | 1.000 | 0.138 | 0.126 | 0.924 |
| | P | −0.001 | 0.031 | 0.028 | 0.918 | −0.990 | 0.133 | 0.125 | 0.936 | 1.000 | 0.139 | 0.126 | 0.930 |
| | N | −0.001 | 0.026 | 0.024 | 0.904 | −0.990 | 0.133 | 0.125 | 0.940 | 0.999 | 0.138 | 0.126 | 0.926 |
| 0.1 | B | 0.002 | 0.028 | 0.027 | 0.920 | −0.890 | 0.139 | 0.126 | 0.910 | 1.003 | 0.129 | 0.126 | 0.936 |
| | P | 0.003 | 0.030 | 0.028 | 0.938 | −0.890 | 0.139 | 0.126 | 0.912 | 1.003 | 0.129 | 0.126 | 0.938 |
| | N | 0.003 | 0.025 | 0.024 | 0.938 | −0.890 | 0.140 | 0.126 | 0.912 | 1.004 | 0.129 | 0.126 | 0.938 |
| 0.2 | B | 0.000 | 0.029 | 0.027 | 0.936 | −0.802 | 0.144 | 0.126 | 0.894 | 0.991 | 0.138 | 0.126 | 0.932 |
| | P | −0.000 | 0.030 | 0.028 | 0.934 | −0.802 | 0.145 | 0.126 | 0.898 | 0.991 | 0.138 | 0.126 | 0.940 |
| | N | −0.001 | 0.027 | 0.024 | 0.912 | −0.801 | 0.145 | 0.126 | 0.896 | 0.992 | 0.138 | 0.125 | 0.934 |
| 0.25 | B | −0.001 | 0.029 | 0.027 | 0.930 | −0.749 | 0.128 | 0.127 | 0.936 | 0.990 | 0.138 | 0.127 | 0.940 |
| | P | −0.000 | 0.031 | 0.028 | 0.928 | −0.748 | 0.129 | 0.127 | 0.938 | 0.990 | 0.139 | 0.127 | 0.938 |
| | N | −0.000 | 0.024 | 0.024 | 0.948 | −0.749 | 0.128 | 0.126 | 0.938 | 0.990 | 0.138 | 0.126 | 0.940 |
| 0.5 | B | −0.002 | 0.029 | 0.027 | 0.926 | −0.513 | 0.143 | 0.126 | 0.918 | 1.000 | 0.131 | 0.126 | 0.936 |
| | P | −0.002 | 0.031 | 0.028 | 0.928 | −0.513 | 0.143 | 0.126 | 0.920 | 1.001 | 0.131 | 0.126 | 0.936 |
| | N | −0.001 | 0.026 | 0.024 | 0.926 | −0.513 | 0.143 | 0.126 | 0.918 | 1.001 | 0.131 | 0.126 | 0.936 |
| 0.7 | B | 0.000 | 0.029 | 0.027 | 0.936 | −0.299 | 0.140 | 0.127 | 0.916 | 0.996 | 0.138 | 0.127 | 0.924 |
| | P | 0.001 | 0.029 | 0.028 | 0.930 | −0.298 | 0.140 | 0.127 | 0.920 | 0.997 | 0.138 | 0.127 | 0.926 |
| | N | 0.001 | 0.025 | 0.024 | 0.934 | −0.299 | 0.140 | 0.126 | 0.914 | 0.996 | 0.138 | 0.126 | 0.926 |
| 1 | B | 0.001 | 0.030 | 0.027 | 0.934 | 0.002 | 0.137 | 0.127 | 0.942 | 1.008 | 0.144 | 0.127 | 0.908 |
| | P | 0.001 | 0.031 | 0.028 | 0.934 | 0.002 | 0.137 | 0.127 | 0.938 | 1.008 | 0.145 | 0.127 | 0.906 |
| | N | 0.001 | 0.026 | 0.024 | 0.928 | 0.002 | 0.138 | 0.127 | 0.938 | 1.007 | 0.144 | 0.127 | 0.908 |

Note: "Est" is the simulation mean; "SE" is the mean of the estimated standard error; "SD" is the mean of the estimated standard deviation; and "COV" is the coverage probability of a nominal 95% confidence interval. The methods used are "B" for the benchmark method, "P" for the proposed method, and "N" for the naive method.

Table 2
*Results of simulation study for scenario* (ii)

| $\varrho$ | | $\beta_1$ | | | | $\beta_2$ | | | | $\beta_3$ | | | |
|---|---|---|---|---|---|---|---|---|---|---|---|---|---|
| | | Est. | SE | SD | COV | Est. | SE | SD | COV | Est. | SE | SD | COV |
| 0 | B | 0.000 | 0.034 | 0.033 | 0.920 | −0.795 | 0.154 | 0.153 | 0.948 | 0.995 | 0.160 | 0.154 | 0.946 |
| | P | 0.001 | 0.036 | 0.035 | 0.922 | −0.795 | 0.154 | 0.154 | 0.950 | 0.995 | 0.159 | 0.154 | 0.950 |
| | N | 0.000 | 0.030 | 0.029 | 0.928 | −0.794 | 0.154 | 0.153 | 0.948 | 0.994 | 0.160 | 0.154 | 0.950 |
| 0.05 | B | 0.002 | 0.028 | 0.027 | 0.920 | −0.790 | 0.139 | 0.125 | 0.910 | 1.003 | 0.129 | 0.126 | 0.936 |
| | P | 0.003 | 0.030 | 0.028 | 0.938 | −0.790 | 0.139 | 0.125 | 0.908 | 1.004 | 0.129 | 0.126 | 0.938 |
| | N | 0.003 | 0.025 | 0.024 | 0.938 | −0.790 | 0.140 | 0.125 | 0.908 | 1.004 | 0.129 | 0.126 | 0.938 |
| 0.1 | B | 0.000 | 0.029 | 0.027 | 0.936 | −0.802 | 0.144 | 0.126 | 0.894 | 0.991 | 0.138 | 0.125 | 0.928 |
| | P | −0.000 | 0.030 | 0.028 | 0.936 | −0.802 | 0.144 | 0.126 | 0.898 | 0.991 | 0.138 | 0.126 | 0.938 |
| | N | −0.001 | 0.027 | 0.024 | 0.916 | −0.801 | 0.145 | 0.125 | 0.896 | 0.992 | 0.138 | 0.125 | 0.932 |
| 0.15 | B | −0.001 | 0.029 | 0.027 | 0.932 | −0.799 | 0.128 | 0.126 | 0.936 | 0.990 | 0.138 | 0.126 | 0.938 |
| | P | −0.000 | 0.031 | 0.028 | 0.930 | −0.798 | 0.128 | 0.126 | 0.938 | 0.990 | 0.138 | 0.126 | 0.938 |
| | N | −0.000 | 0.024 | 0.024 | 0.950 | −0.799 | 0.128 | 0.126 | 0.938 | 0.990 | 0.138 | 0.126 | 0.936 |
| 0.2 | B | −0.002 | 0.029 | 0.027 | 0.926 | −0.813 | 0.143 | 0.126 | 0.918 | 1.001 | 0.131 | 0.126 | 0.934 |
| | P | −0.002 | 0.031 | 0.028 | 0.932 | −0.813 | 0.143 | 0.126 | 0.918 | 1.001 | 0.131 | 0.126 | 0.934 |
| | N | −0.001 | 0.026 | 0.024 | 0.924 | −0.813 | 0.143 | 0.126 | 0.916 | 1.001 | 0.131 | 0.126 | 0.934 |
| 0.5 | B | 0.000 | 0.029 | 0.027 | 0.936 | −0.799 | 0.140 | 0.126 | 0.916 | 0.996 | 0.138 | 0.126 | 0.924 |
| | P | 0.001 | 0.029 | 0.028 | 0.930 | −0.798 | 0.140 | 0.126 | 0.922 | 0.997 | 0.138 | 0.127 | 0.926 |
| | N | 0.001 | 0.025 | 0.024 | 0.934 | −0.799 | 0.140 | 0.126 | 0.914 | 0.996 | 0.138 | 0.126 | 0.926 |
| 0.7 | B | 0.001 | 0.030 | 0.027 | 0.932 | −0.798 | 0.137 | 0.127 | 0.942 | 1.008 | 0.144 | 0.127 | 0.906 |
| | P | 0.001 | 0.031 | 0.028 | 0.934 | −0.798 | 0.137 | 0.127 | 0.938 | 1.008 | 0.145 | 0.127 | 0.906 |
| | N | 0.001 | 0.026 | 0.024 | 0.930 | −0.798 | 0.138 | 0.126 | 0.938 | 1.007 | 0.144 | 0.126 | 0.908 |

TABLE 3
*Results of simulation study for scenario* (iii)

| $\varrho$ | | $\beta_1$ | | | | $\beta_2$ | | | | $\beta_3$ | | | |
|---|---|---|---|---|---|---|---|---|---|---|---|---|---|
| | | Est. | SE | SD | COV | Est. | SE | SD | COV | Est. | SE | SD | COV |
| 0 | B | 0.200 | 0.034 | 0.035 | 0.936 | −0.995 | 0.154 | 0.162 | 0.956 | 0.995 | 0.160 | 0.163 | 0.964 |
| | P | 0.195 | 0.038 | 0.038 | 0.920 | −0.995 | 0.158 | 0.167 | 0.958 | 0.994 | 0.163 | 0.168 | 0.958 |
| | N | 0.156 | 0.031 | 0.030 | 0.684 | −0.995 | 0.158 | 0.160 | 0.952 | 0.995 | 0.165 | 0.160 | 0.946 |
| 0.05 | B | 0.202 | 0.028 | 0.029 | 0.948 | −0.990 | 0.139 | 0.138 | 0.938 | 1.003 | 0.129 | 0.138 | 0.950 |
| | P | 0.197 | 0.032 | 0.032 | 0.944 | −0.994 | 0.144 | 0.144 | 0.938 | 1.004 | 0.139 | 0.144 | 0.948 |
| | N | 0.159 | 0.026 | 0.025 | 0.602 | −0.991 | 0.149 | 0.133 | 0.910 | 1.004 | 0.140 | 0.133 | 0.926 |
| 0.1 | B | 0.200 | 0.029 | 0.029 | 0.950 | −1.002 | 0.144 | 0.138 | 0.924 | 0.991 | 0.138 | 0.138 | 0.956 |
| | P | 0.194 | 0.032 | 0.033 | 0.948 | −1.005 | 0.151 | 0.144 | 0.922 | 0.991 | 0.147 | 0.144 | 0.954 |
| | N | 0.155 | 0.028 | 0.025 | 0.560 | −1.004 | 0.153 | 0.133 | 0.904 | 0.991 | 0.148 | 0.133 | 0.912 |
| 0.15 | B | 0.199 | 0.029 | 0.029 | 0.950 | −0.999 | 0.128 | 0.138 | 0.960 | 0.990 | 0.138 | 0.138 | 0.960 |
| | P | 0.194 | 0.033 | 0.032 | 0.938 | −0.998 | 0.135 | 0.144 | 0.958 | 0.986 | 0.144 | 0.144 | 0.956 |
| | N | 0.155 | 0.025 | 0.025 | 0.542 | −0.997 | 0.138 | 0.133 | 0.938 | 0.989 | 0.145 | 0.133 | 0.948 |
| 0.2 | B | 0.198 | 0.029 | 0.029 | 0.948 | −1.013 | 0.143 | 0.138 | 0.942 | 1.001 | 0.131 | 0.138 | 0.958 |
| | P | 0.193 | 0.033 | 0.032 | 0.936 | −1.012 | 0.148 | 0.144 | 0.938 | 0.997 | 0.136 | 0.144 | 0.954 |
| | N | 0.155 | 0.027 | 0.025 | 0.536 | −1.016 | 0.154 | 0.133 | 0.920 | 0.998 | 0.141 | 0.133 | 0.932 |
| 0.5 | B | 0.200 | 0.029 | 0.029 | 0.956 | −0.999 | 0.140 | 0.138 | 0.954 | 0.996 | 0.138 | 0.138 | 0.944 |
| | P | 0.195 | 0.032 | 0.032 | 0.952 | −1.000 | 0.147 | 0.144 | 0.954 | 0.993 | 0.144 | 0.144 | 0.956 |
| | N | 0.157 | 0.026 | 0.025 | 0.582 | −1.000 | 0.153 | 0.133 | 0.898 | 0.996 | 0.145 | 0.133 | 0.920 |
| 0.7 | B | 0.201 | 0.030 | 0.029 | 0.952 | −0.998 | 0.137 | 0.139 | 0.958 | 1.008 | 0.144 | 0.139 | 0.938 |
| | P | 0.196 | 0.033 | 0.033 | 0.946 | −0.997 | 0.143 | 0.145 | 0.962 | 1.008 | 0.146 | 0.145 | 0.956 |
| | N | 0.157 | 0.028 | 0.025 | 0.594 | −1.000 | 0.146 | 0.134 | 0.932 | 1.006 | 0.151 | 0.134 | 0.912 |

Table 4
*Results of simulation study for scenario* (iv)

| $\varrho$ | | $\beta_1$ | | | | $\beta_2$ | | | | $\beta_3$ | | | |
|---|---|---|---|---|---|---|---|---|---|---|---|---|---|
| | | Est. | SE | SD | COV | Est. | SE | SD | COV | Est. | SE | SD | COV |
| 0.30 | B | 0.194 | 0.038 | 0.046 | 0.970 | −0.987 | 0.145 | 0.159 | 0.980 | 1.005 | 0.149 | 0.156 | 0.980 |
| | P | 0.173 | 0.040 | 0.042 | 0.850 | −0.986 | 0.150 | 0.152 | 0.970 | 1.005 | 0.150 | 0.149 | 0.970 |
| | N | 0.073 | 0.025 | 0.025 | 0.000 | −0.976 | 0.159 | 0.137 | 0.920 | 0.996 | 0.153 | 0.135 | 0.920 |
| 0.40 | B | 0.199 | 0.043 | 0.044 | 0.950 | −1.002 | 0.123 | 0.147 | 0.970 | 1.002 | 0.131 | 0.147 | 0.980 |
| | P | 0.185 | 0.045 | 0.041 | 0.890 | −1.003 | 0.124 | 0.144 | 0.970 | 0.999 | 0.130 | 0.144 | 0.980 |
| | N | 0.096 | 0.029 | 0.028 | 0.060 | −1.002 | 0.127 | 0.134 | 0.970 | 1.002 | 0.135 | 0.134 | 0.960 |
| 0.50 | B | 0.199 | 0.043 | 0.042 | 0.960 | −0.981 | 0.134 | 0.142 | 0.960 | 1.020 | 0.122 | 0.143 | 0.970 |
| | P | 0.190 | 0.044 | 0.040 | 0.920 | −0.981 | 0.133 | 0.139 | 0.950 | 1.019 | 0.127 | 0.141 | 0.950 |
| | N | 0.116 | 0.033 | 0.030 | 0.200 | −0.988 | 0.137 | 0.133 | 0.930 | 1.020 | 0.132 | 0.135 | 0.930 |
| 0.60 | B | 0.194 | 0.035 | 0.040 | 0.970 | −0.993 | 0.136 | 0.141 | 0.950 | 1.025 | 0.137 | 0.138 | 0.930 |
| | P | 0.192 | 0.038 | 0.040 | 0.950 | −0.994 | 0.140 | 0.140 | 0.950 | 1.025 | 0.138 | 0.137 | 0.910 |
| | N | 0.131 | 0.028 | 0.032 | 0.450 | −0.998 | 0.152 | 0.137 | 0.910 | 1.020 | 0.138 | 0.134 | 0.940 |
| 0.70 | B | 0.198 | 0.038 | 0.039 | 0.960 | −1.018 | 0.137 | 0.133 | 0.970 | 1.004 | 0.140 | 0.131 | 0.930 |
| | P | 0.194 | 0.040 | 0.038 | 0.950 | −1.017 | 0.138 | 0.132 | 0.960 | 1.004 | 0.142 | 0.131 | 0.930 |
| | N | 0.152 | 0.038 | 0.033 | 0.660 | −1.021 | 0.142 | 0.130 | 0.920 | 1.004 | 0.144 | 0.128 | 0.920 |
| 0.80 | B | 0.203 | 0.036 | 0.038 | 0.950 | −1.001 | 0.142 | 0.132 | 0.930 | 1.005 | 0.136 | 0.132 | 0.960 |
| | P | 0.203 | 0.038 | 0.038 | 0.950 | −1.002 | 0.143 | 0.131 | 0.940 | 1.005 | 0.135 | 0.132 | 0.960 |
| | N | 0.172 | 0.035 | 0.035 | 0.870 | −1.002 | 0.147 | 0.131 | 0.920 | 1.000 | 0.136 | 0.131 | 0.930 |

4.2. *Performance of the proposed tests.* We now explore the numerical performance of the proposed tests. First, we want to test a hypothesis of the parametric component of form:

(4.5) $$H_0 : A\boldsymbol{\beta} = 0 \quad \text{VS} \quad H_1 : A\boldsymbol{\beta} = c,$$

where $A = (1,1,1)^{\mathrm{T}}$, $c$ is a value from the set $\{0, 0.1, 0.2, \ldots, 0.7, 1\}$, $\boldsymbol{\beta} = (0.2, c-1.2, 1)^{\mathrm{T}}$ and $\alpha_1(\cdot)$ and $\alpha_1(\cdot)$ are the same as those in scenario (i). The same models and error distribution as in Section 4.1 are used.

The power to detect $H_1$ was calculated by using the critical values from the chi-squared approximation and the wild bootstrap approximation. To compare test performances, the powers of the tests based on the benchmark estimator, the proposed estimator and the naive estimator are presented. In implementing the wild bootstrap method, we generated 500 bootstrap samples from the model

$$\begin{cases} Y_i^* = \hat{\beta}_1 \hat{\xi}_i + \hat{\beta}_2 W_{1i} + \hat{\beta}_3 W_{2i} + \hat{\alpha}_1(U_i) Z_{1i} + \hat{\alpha}_2(U_i) Z_{2i} + \varepsilon_i^*, \\ \hat{\xi}_i = \hat{\xi}(V_i), \end{cases}$$

where, $\varepsilon_i^*$ is a wild bootstrap residual; that is, $\varepsilon_i^* = \tau_i \hat{\varepsilon}_i$, with $\hat{\varepsilon}_i = Y_i - \{\hat{\beta}_1 \hat{\xi}_i + \hat{\beta}_2 W_{1i} + \hat{\beta}_3 W_{2i} + \hat{\alpha}_1(U_i) Z_{1i} + \hat{\alpha}_2(U_i) Z_{2i}\}$, $\tau_i = -(\sqrt{5}-1)/2$ with probability $(\sqrt{5}+1)/(2\sqrt{5})$ and $\tau_i = (\sqrt{5}+1)/2$ with $1 - (\sqrt{5}+1)/(2\sqrt{5})$. Using this bootstrap sample $(Y_i^*, \hat{\xi}_i, \mathbf{W}_i, \mathbf{Z}_i, U_i)$, we can calculate the $T_n^*$ and $W_n^*$, and get the 95 percentiles as the critical values for the proposed tests at the significance level 0.05.

The power of $T_n$ associated to scenario (iii) is presented in Table 5 for $\beta_1 = 0.2$. Note that the power is actually the empirical level when $c = 0$. All empirical levels close nominal level 0.05 and the empirical level based on the wild bootstrap procedure are consistently smaller than those based on the $\mathcal{X}^2$ approximation and are closer to the nominal level. These facts apply for $\beta_1 = 0$ (not listed here). As $c$ increases to 0.7, the powers of two tests based on $\mathcal{X}^2$ approximation are greater than 0.92. Similar conclusions can be drawn for the Wald test, whose simulation results are also given in Table 5.

We further study the numerical performance of the test by checking the nonparametric component. We consider the following hypothesis:

(4.6) $$H_0 : \alpha_1(u) = m \quad \text{VS} \quad \alpha_1(u) = \alpha_1(u, \varrho) \text{ given by } (4.3).$$

The simulation results obtained by using the wild bootstrap approximation method to choose critical value are shown in Table 6. When $\varrho = 0$, the results are the empirical levels, which are close to the nominal level. The power is greater than 0.99 when $\varrho = 0.5$. Table 6 also indicates that the power is a monotone increasing function of $\varrho$.

TABLE 5
*Empirical power of profile least-square ratio test $T_n$ and the Wald test $W_n$ at level 0.05 for hypothesis (4.5). Data were generated from models (4.1) with $\boldsymbol{\beta} = (0.2, c-1.2, 1)^{\mathrm{T}}$ and $c \in \{0, 0.1, 0.2, 0.25, 0.5, 0.7, 1\}$ and $\alpha_1(u)$ and $\alpha_2(u)$ given by (4.2) and (4.4), respectively. The methods used are "Asm" for the asymptotic version, and "Boot" for the bootstrap version*

| | | $T_n$ | | | Wald | | |
|---|---|---|---|---|---|---|---|
| $c$ | | B | P | N | B | P | N |
| 0 | Aym | 0.060 | 0.070 | 0.080 | 0.050 | 0.050 | 0.080 |
| | Boot | 0.050 | 0.060 | 0.060 | 0.060 | 0.060 | 0.060 |
| 0.10 | Aym | 0.150 | 0.140 | 0.150 | 0.130 | 0.130 | 0.150 |
| | Boot | 0.130 | 0.100 | 0.080 | 0.130 | 0.120 | 0.080 |
| 0.20 | Aym | 0.190 | 0.220 | 0.120 | 0.150 | 0.150 | 0.120 |
| | Boot | 0.170 | 0.160 | 0.080 | 0.190 | 0.180 | 0.080 |
| 0.25 | Aym | 0.350 | 0.340 | 0.240 | 0.320 | 0.310 | 0.240 |
| | Boot | 0.290 | 0.280 | 0.180 | 0.310 | 0.300 | 0.180 |
| 0.50 | Aym | 0.740 | 0.710 | 0.530 | 0.670 | 0.660 | 0.530 |
| | Boot | 0.700 | 0.630 | 0.500 | 0.720 | 0.630 | 0.500 |
| 0.70 | Aym | 0.940 | 0.940 | 0.870 | 0.930 | 0.920 | 0.870 |
| | Boot | 0.920 | 0.890 | 0.800 | 0.930 | 0.890 | 0.800 |
| 1.00 | Aym | 1.000 | 1.000 | 1.000 | 0.990 | 0.990 | 1.000 |
| | Boot | 0.990 | 0.990 | 0.960 | 0.990 | 0.990 | 0.960 |

TABLE 6
*Empirical power of level 0.05 for hypothesis (4.6) using the wild bootstrap procedure. Data were generated from (4.1) and (4.3) with $\boldsymbol{\beta} = (0.2, -1, 1)^{\mathrm{T}}$ and $\varrho \in \{0, 0.5, 0.10, 0.15, 0.5, 0.7\}$*

| $\varrho$ | B | P | N |
|---|---|---|---|
| 0 | 0.060 | 0.050 | 0.080 |
| 0.05 | 0.110 | 0.140 | 0.160 |
| 0.10 | 0.240 | 0.260 | 0.250 |
| 0.15 | 0.410 | 0.360 | 0.360 |
| 0.20 | 0.520 | 0.510 | 0.500 |
| 0.50 | 0.990 | 0.990 | 1.000 |
| 0.70 | 1.000 | 1.000 | 1.000 |

4.3. *Real data example.* To illustrate the proposed estimation method, we consider a dataset from a Duchenne Muscular Dystrophy (DMD) study. See Andrews and Herzberg (1985) for a detailed discussion on the dataset. The dataset contains 209 observations corresponding to blood samples on 192 patients (17 patients have two samples) collected from a project to develop a screening program for female relatives of boys with DMD. The program's goal was to inform a woman of her chances of being a carrier based on serum markers as well as her family pedigree. Another question of interest is whether age should be taken into account in the analysis. Enzyme levels were measured in known carriers (75 samples) and in a group of noncarriers (134 samples). The serum marker creatine kinase (ck) is inexpensive to obtain, while the marker lactate dehydrogenase (ld) is very expensive to obtain. It is of interest to predict the value ld by using the level of ck, carrier status and age of patient.

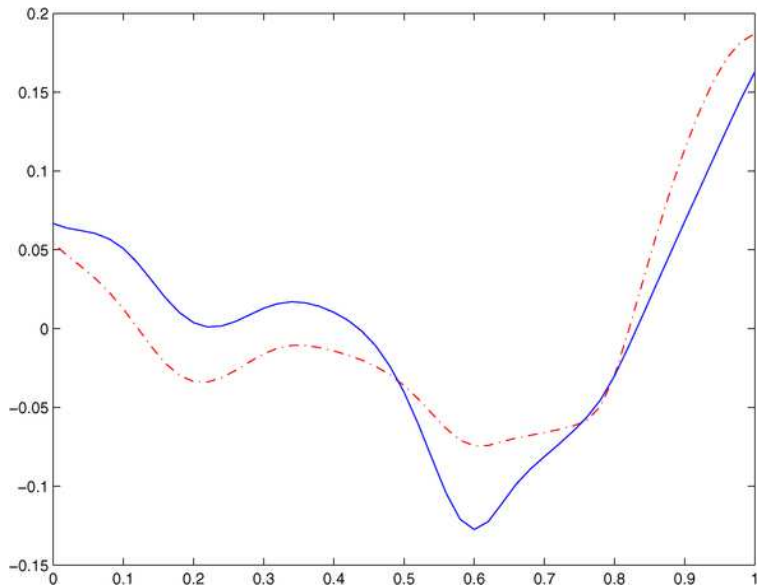

Fig. 1. *Estimated curves of the nonparametric function for the DMD study. The solid, dotted lines were obtained using the naive and proposed methods, respectively.*

We consider the following model: $Y = \beta_0 + \beta_1 Z_1 + \beta_2 Z_2 + g(U)$, where $Z_1 = $ ck is measured with errors and $Z_2 = $ carrier status is exactly measured, $U$ is age and $Y$ denotes the observed level of lactate dehydrogenase. We justify the measurement error of $Z_1$ by regressing $Z_1$ on $U$. The estimates and associated standard errors based on the naive and proposed methods are as follows: $\widehat{\beta}_{0,\text{naive}} = 4.6057(0.113)$, $\widehat{\beta}_{1,\text{naive}} = 0.1509(0.027)$ and $\widehat{\beta}_{2,\text{naive}} = 0.2269(0.055)$; $\widehat{\beta}_{0,n} = 4.4296(0.329)$, $\widehat{\beta}_{1,n} = 0.1775(0.042)$ and $\widehat{\beta}_{2,n} = 0.3702(0.050)$. The estimated curves of the nonparametric function $g(u)$ are provided in Figure 1. Accounting for measurement errors, the estimate of $\beta_1$ increases about 17.2%, and the associated standard error also increases 55%. The estimate of $\beta_2$ also increases when measurement errors are taken into account. The patterns of the nonparametric curve are similar, and show a slight difference.

**5. Discussion.** We developed estimation and inference procedures for the SVCPLM when parts of the parametric components are unobserved. The procedures are derived by incorporating ancillary information to calibrate the mismeasured variables and by applying the profile least-square-based principle.

In some cases we may not have an auxiliary variable $\boldsymbol{\eta}$, but we can observe two or more independent replicates of $\mathbf{V}$. For instance, when two measurements $\mathbf{V}_1$ and $\mathbf{V}_2$, which satisfy that $\mathbf{V}_1 = \boldsymbol{\xi} + u_1$ and $\mathbf{V}_2 = \boldsymbol{\xi} + u_2$, and $E(u_1|\mathbf{V}_2) = 0$ and $E(u_2|\mathbf{V}_1) = 0$, are available, we can estimate $\boldsymbol{\xi}$ by

$$\hat{\boldsymbol{\xi}}(v) = \frac{\sum_{i=1}^n \{\mathbf{V}_{i1} K_h(\mathbf{V}_{i2} - v) + \mathbf{V}_{i1} K_h(\mathbf{V}_{i1} - v)\}}{\sum_{i=1}^n \{K_h(\mathbf{V}_{i2} - v) + K_h(\mathbf{V}_{i1} - v)\}},$$

because $E(\mathbf{V}_1|\mathbf{V}_2 = v) = E(\mathbf{V}_2|\mathbf{V}_1 = v) = E(\boldsymbol{\xi}|\mathbf{V} = v)$. The proposed procedure applies to this situation as well, and similar results to those presented in this paper can be obtained for the resulting estimator.

It is of interest to extend the proposed methodology to a more general semiparametric model: $E(Y|\mathbf{Z}, \mathbf{X}, U) = G\{\boldsymbol{\Theta}^{\mathrm{T}}\mathbf{Z} + \boldsymbol{\alpha}^{\mathrm{T}}(U)\mathbf{X}\}$, where $G(\cdot)$ is a link function. The study of this model with mismeasured components of $\mathbf{Z}$ needs further investigation and is beyond the scope of this paper.

## APPENDIX

In this Appendix, we list assumptions and outline proofs of the main results. The following technical assumptions are imposed:

### A.1. Assumptions.

1. The random variable $U$ has a bounded support $\mathcal{U}$. Its density function $f_u(\cdot)$ is Lipschitz continuous and bounded away from 0 on its support. The density function of random variable $\mathbf{V}$, $f_v(v)$, is continuously differentiable and bounded away from 0 and infinite on its finite support $\mathcal{V}$. $\{\alpha_i(u), i = 1, 2, \ldots, q\}$ have a continuous second derivative.
2. The $q \times q$ matrix $E(\mathbf{ZZ}^{\mathrm{T}}|U)$ is nonsingular for each $U \in \mathcal{U}$. All elements of the matrices $E(\mathbf{ZZ}^{\mathrm{T}}|U)$, $E(\mathbf{ZZ}^{\mathrm{T}}|U)^{-1}$ and $E(\mathbf{ZX}^{\mathrm{T}}|U)$ are Lipschitz continuous.
3. The kernel functions $K(\cdot)$ and $L(\cdot)$ are density functions with compact support $[-1, 1]$.
4. There is an $s > 2$ such that $E\|\mathbf{Z}\|^{2s} < \infty$ and $E\|\mathbf{X}\|^{2s} < \infty$ and for some $\delta < 2 - s^{-1}$ such that $n^{2\delta-1}h \to \infty$, $n^{2\delta-1}b_k \to \infty$ and $nhb_k^{(2r+2)} \to 0$, $k = 1, 2, \ldots, p_1$, where $b_k$ is the bandwidth parameter in the polynomial estimator $\hat{\xi}_k(\cdot)$ of $\xi_k(\cdot)$.
5. $nh^8 \to 0$ and $nh^2/(\log n)^2 \to \infty$.

**A.2. Preliminary lemmas.** Write $c_{n1} = (\frac{\log h}{nh})^{1/2} + h^2$, $c_{n2} = (\frac{\log b}{nb})^{1/2} + b^{r+1}$, $c_n = c_{n1} + c_{n2}$.

LEMMA A.1. *Suppose that $(\mathbf{Z}_i, \mathbf{X}_i, U_i), i = 1, 2, \ldots, n$ are an i.i.d. random vector. $E|g(\mathbf{X}, \mathbf{Z}, U)| < \infty$ and $E[g(\cdot, \cdot, u)|U = u]$ have a continuous second derivative on $u$. Further assume that $E(|g(\mathbf{X}, \mathbf{Z}, U)|^s | \mathbf{Z} = \mathbf{z}, \mathbf{X} = \mathbf{x}) < \infty$. Let $K$ be a bounded positive function with a bounded support satisfying the Lipschitz condition. Given that $n^{2\delta-1}h \to \infty$ for some $\delta < 1 - s^{-1}$, then we have*

$$\sup_{u \in \mathcal{U}} \left| \frac{1}{n} \sum_{i=1}^n K_h(U_i - u) \left(\frac{U_i - u}{h}\right)^k g(X_i, \mathbf{Z}_i, U_i) - f(u) E\{g(\mathbf{X}, \mathbf{Z}, u) | U = u\} \mu_k \right|$$
$$= O(c_{n1}) \quad a.s.$$

*Furthermore, assume that $E[\varepsilon_i | \mathbf{Z}_i, \mathbf{X}_i, U_i] = 0$, $E[|\varepsilon_i|^s | \mathbf{Z}_i, \mathbf{X}_i, U_i)] < \infty$, then*

$$\sup_{u \in \mathcal{U}} \left| \frac{1}{n} \sum_{i=1}^n K_h(U_i - u) g(\mathbf{X}_i, Z_i, U_i) \varepsilon_i \right| = O(c_{n1}) \quad a.s.$$

PROOF. The first result follows an argument similar to that of Lemma A.2 of Fan and Huang (2005). The second result follows the first result and an argument similar to Xia and Li (1999). □

LEMMA A.2. *Suppose that $E[g(\mathbf{Z}, \mathbf{X}, u)|U = u]$ has a continuous second derivative on $u$ and $E|g(\mathbf{X}, \mathbf{Z}, U)|^s < \infty$. Under Assumptions 1–5, we have*

$$\sup_{u \in \mathcal{U}} \left| \frac{1}{n} \sum_{i=1}^n K_h(U_i - u) \left(\frac{U_i - u}{h}\right)^k g(\mathbf{X}_i, \mathbf{Z}_i, U_i) \hat{\boldsymbol{\xi}}_i^{\mathrm{T}} \right.$$

$$\left. - f(u) E\{g(\mathbf{X}, \mathbf{Z}, u) \boldsymbol{\xi}^{\mathrm{T}} | U = u\} \mu_k \right| = O(c_n) \quad a.s.$$

*and*

$$\sup_{u \in \mathcal{U}} \left| \frac{1}{n} \sum_{i=1}^n K_h(U_i - u) g(\mathbf{X}_i, \mathbf{Z}_i, U_i) h(\hat{\boldsymbol{\xi}}_i) \varepsilon_i \right| = O(c_n),$$

*where $h(\cdot)$ is a twice continuous differentiable function.*

PROOF. Note that $\frac{1}{n} \sum_{i=1}^n K_h(U_i - u)(\frac{U_i - u}{h})^k g(\mathbf{X}_i, \mathbf{Z}_i, U_i) \hat{\boldsymbol{\xi}}_i^{\mathrm{T}}$ can be decomposed as

$$\frac{1}{n} \sum_{i=1}^n K_h(U_i - u) \left(\frac{U_i - u}{h}\right)^k g(\mathbf{X}_i, \mathbf{Z}_i, U_i) \boldsymbol{\xi}_i^{\mathrm{T}}$$

$$+ \frac{1}{n} \sum_{i=1}^n K_h(U_i - u) \left(\frac{U_i - u}{h}\right)^k g(\mathbf{X}_i, \mathbf{Z}_i, U_i) (\hat{\boldsymbol{\xi}}_i - \boldsymbol{\xi}_i)^{\mathrm{T}}.$$

By Lemma A.1, the first term equals $f_u(u)E\{g(\mathbf{X},\mathbf{Z},u)\boldsymbol{\xi}|U=u\}\mu_k + O(c_{n1})$ uniformly on $u \in \mathcal{U}$ in probability. Recalling the asymptotic expression given in (2.1) and using Lemma A.1, one can show that the second term is $O(c_{n2})$. This completes the proof of Lemma 2. $\square$

LEMMA A.3. $g(\cdot,\cdot,u)$ has a continuous second derivative on $u$ and $E|g(\mathbf{X}, \mathbf{Z}, U)| < \infty$. Under Assumptions 1–5, $n^{-1}\sum_{i=1}^n (\hat{\mathbf{Z}}_i - \mathbf{Z}_i)\hat{\mathbf{Z}}_i^l g(\mathbf{X}_i, \mathbf{Z}_i, U_i)$ is of order $O(c_n)$ a.s., where $l = 0, 1$.

PROOF. The proof follows from (2.1) and arguments similar to Lemma A.2. $\square$

LEMMA A.4. Under Assumptions 1–5, we have

$$(\tilde{\mathbf{Z}}^T\tilde{\mathbf{Z}})^{-1}\tilde{\mathbf{Z}}^T(\mathbf{I} - \mathbf{S})\mathbf{Z}$$
$$= \boldsymbol{\Sigma}^{-1}\left(\frac{\zeta_1 S_u^{-1} c_p b^{r+1}}{n(r+1)!}\sum_{i=1}^n \psi(\mathbf{Z}_i, \mathbf{X}_i, U_i)[\{\boldsymbol{\xi}^{(r+1)}(\mathbf{V}_i)\}^T, \mathbf{0}]\right.$$
$$\left. + \frac{1}{n^2}\sum_{i=1}^n\sum_{j=1}^n \frac{1}{f_v(\mathbf{V}_i)}\psi(\mathbf{Z}_i, \mathbf{X}_i, U_i) L_b(\mathbf{V}_j - \mathbf{V}_i)(\mathbf{e}_j^T, \mathbf{0})\right)\{1 + o(1)\}$$

in probability.

PROOF. We first prove that

(A.1) $$\frac{1}{n}\tilde{\mathbf{Z}}^T\tilde{\mathbf{Z}} \to \boldsymbol{\Sigma}.$$

A direct calculation yields

(A.2) $$D_u^T W_u D_u = nf_u(U)\Gamma(U) \otimes \begin{pmatrix} 1 & \mu_1 \\ \mu_1 & \mu_2 \end{pmatrix}\{1 + O_P(c_{n1})\}.$$

On the other hand, Lemma A.3 implies

(A.3) $$D_u^T W_u \hat{\mathbf{Z}} = nf_u(U)\Phi(U) \otimes (1,\mu_1)^T\{1 + O_P(c_n)\}.$$

A combination of (A.2) and (A.3) implies

(A.4) $(\mathbf{X}^T, 0)(D_u^T W_u D_u)^{-1} D_u^T W_u \hat{\mathbf{Z}} = \mathbf{X}^T \Gamma^{-1}(U)\Phi(U)\{1 + O_P(c_n)\}$

and then

(A.5) $\tilde{\mathbf{Z}}_i = \hat{\mathbf{Z}}_i - \Phi^T(U_i)\Gamma^{-1}(U_i)\mathbf{X}_i\{1 + O_P(c_n)\}, \qquad i = 1, 2, \ldots, n.$

It follows from these arguments that $n^{-1}\tilde{\mathbf{Z}}^T\tilde{\mathbf{Z}} = \frac{1}{n}\sum_{i=1}^n \{\psi(\mathbf{Z}_i, \mathbf{X}_i, U_i)\}^{\otimes 2}\{1 + O_P(c_n)\}$, and (A.1) follows.

Note that $\tilde{\mathbf{Z}}^{\mathrm{T}}(\mathbf{I}-\mathbf{S})(\mathbf{Z}-\hat{\mathbf{Z}}) = \mathbf{Z}^{\mathrm{T}}(\mathbf{I}-\mathbf{S})^{\mathrm{T}}(\mathbf{I}-\mathbf{S})(\mathbf{Z}-\hat{\mathbf{Z}}) - (\mathbf{Z}-\hat{\mathbf{Z}})^{\mathrm{T}}(\mathbf{I}-\mathbf{S})^{\mathrm{T}}(\mathbf{I}-\mathbf{S})(\mathbf{Z}-\hat{\mathbf{Z}}) \overset{\mathrm{def}}{=} J_1 - J_2$. The second term, $J_2$, is $O_{\mathrm{P}}(c_n^2)$ by Lemma A.3. Write $\tilde{\mathbf{Z}}_* = (\mathbf{I}-\mathbf{S})\mathbf{Z}$. We have $J_1 = \tilde{\mathbf{Z}}_*^{\mathrm{T}}(\mathbf{Z}-\hat{\mathbf{Z}}) - \tilde{\mathbf{Z}}_*^{\mathrm{T}}\mathbf{S}(\mathbf{Z}-\hat{\mathbf{Z}})$. It follows from (2.1) that

$$D_u W_u(\mathbf{Z}-\hat{\mathbf{Z}})$$

$$= \frac{\zeta_1 S_u^{-1} c_p b^{r+1}}{(r+1)!} \sum_{i=1}^n K_h(U_i - U)\mathbf{X}_i \boldsymbol{\xi}^{(r+1)}(\mathbf{V}_i) \otimes \begin{pmatrix} 1 & 0 \\ U_i - U & 0 \\ h & \end{pmatrix}$$

$$+ \frac{1}{n^2}\sum_{i=1}^n\sum_{j=1}^n f_v^{-1}(V_i) K_h(U_i - U) L_b(V_j - V_i)\mathbf{X}_i \mathbf{e}_j^{\mathrm{T}} \otimes \begin{pmatrix} 1 & 0 \\ U_i - U & 0 \\ h & \end{pmatrix}$$

$$+ o(b^{r+1} + \log b^{-1}/\sqrt{nb}).$$

By an argument similar to that of (A.5), we derive

$$\tilde{\mathbf{Z}}_*^{\mathrm{T}}\mathbf{S}(\mathbf{Z}-\hat{\mathbf{Z}})$$

$$= \frac{1}{n}\sum_{l=1}^n \tilde{\rho}(\mathbf{Z}_l, \mathbf{X}_l, U_l)$$

$$\times \sum_{i=1}^n \left\{ \frac{\zeta_1 S_u^{-1} c_p b^{r+1}}{(r+1)!} K_h(U_i - U_l)\mathbf{X}_i \boldsymbol{\xi}^{(r+1)}(\mathbf{V}_i) \otimes \begin{pmatrix} 1 & 0 \\ U_i - U_l & 0 \\ h & \end{pmatrix} \right.$$

$$+ \frac{1}{n}\sum_{j=1}^n f_v^{-1}(V_i) K_h(U_i - U_l) L_b(V_j - V_i)\mathbf{X}_i \mathbf{e}_j^{\mathrm{T}}$$

$$\left. \otimes \begin{pmatrix} 1 & 0 \\ U_i - U_l & 0 \\ h & \end{pmatrix} \right\}$$

$$\times \{1 + o_{\mathrm{P}}(1)\},$$

where $\tilde{\rho}(\mathbf{Z}_l, \mathbf{X}_l, U_l)$ can be expressed as

$$\psi(\mathbf{Z}_l, \mathbf{X}_l, U_l)\{1 + O_{\mathrm{P}}(c_n)\}(\mathbf{X}_l^{\mathrm{T}}, \mathbf{0})$$

$$\times \left\{ f_u(U_l)\Gamma(U_l) \otimes \begin{pmatrix} 1 & \mu_1 \\ \mu_1 & \mu_2 \end{pmatrix} \{1 + O_{\mathrm{P}}(c_n)\} \right\}^{-1}$$

$$= \frac{1}{f_u(U_l)}\{\psi(\mathbf{Z}_l, \mathbf{X}_l, U_l)^{\mathrm{T}}\mathbf{X}_l^{\mathrm{T}}, 0\}$$

$$\times \left\{ \Gamma^{-1}(U_l) \otimes \begin{pmatrix} 1 & \mu_1 \\ \mu_1 & \mu_2 \end{pmatrix}^{-1} \right\}\{1 + O_{\mathrm{P}}(c_n)\}$$

$$= \frac{\psi(\mathbf{Z}_l, \mathbf{X}_l, U_l)\mathbf{X}_l^{\mathrm{T}}}{f_u(U_l)(\mu_2 - \mu_1^2)}\Gamma^{-1}(U_l) \otimes (\mu_2, -\mu_1)\{1 + O_{\mathrm{P}}(c_n)\}.$$

Denote by $\rho(\mathbf{Z}_l, \mathbf{X}_l, U_l)$ the main term of the right-hand side of the above formula. Note that $E\{\rho(\mathbf{Z}_l, \mathbf{X}_l, U_l)|U_l\} = 0$. By Lemma 3 of Chen, Choi and Zhou (2005) we have

$$\frac{1}{n^3} \sum_{i=1}^{n} \sum_{j=1}^{n} \sum_{l=1}^{n} K_h(U_i - U_l) L_b(V_j - V_i) \rho(\mathbf{Z}_l, \mathbf{X}_l, U_l) \frac{\mathbf{X}_i \mathbf{e}_j^{\mathrm{T}}}{f_v(V_i)}$$

(A.6)
$$= O_{\mathrm{P}}(c_n n^{-1/2}).$$

Furthermore, we can show in a similar way as that for (A.6), that

$$\frac{\zeta_1 S_u^{-1} c_p b^{2(r+1)}}{n^2 (r+1)!} \sum_{i=1}^{n} \sum_{l=1}^{n} K_h(U_i - U_l) \rho(\mathbf{Z}_l, \mathbf{X}_l, U_l) \mathbf{X}_i \{\boldsymbol{\xi}^{(r+1)}(V_i)\}^{\mathrm{T}} = O_{\mathrm{P}}(c_n^2).$$

These arguments imply that

(A.7) $$n^{-1} \tilde{\mathbf{Z}}_*^{\mathrm{T}} \mathbf{S}(\mathbf{Z} - \hat{\mathbf{Z}}) = O_{\mathrm{P}}(c_n^2).$$

We now deal with the term $\tilde{\mathbf{Z}}_*^{\mathrm{T}}(\mathbf{Z}-\hat{\mathbf{Z}})$. Note that $\tilde{\mathbf{Z}}_*^{\mathrm{T}}(\mathbf{Z}-\hat{\mathbf{Z}})$ equals $\sum_{i=1}^{n} \psi(\mathbf{Z}_i, \mathbf{X}_i, U_i)\{(\boldsymbol{\xi}_i - \hat{\boldsymbol{\xi}}_i)^{\mathrm{T}}, \mathbf{0}\}$, which can be further decomposed as

$$\frac{\zeta_1 S_u^{-1} c_p b^{r+1}}{(r+1)!} \sum_{i=1}^{n} \psi(\mathbf{Z}_i, \mathbf{X}_i, U_i)[\{\boldsymbol{\xi}^{(r+1)}(V_i)\}^{\mathrm{T}}, \mathbf{0}]$$
$$+ \frac{1}{n} \sum_{i=1}^{n} \sum_{j=1}^{n} \frac{1}{f_v(V_i)} \psi(\mathbf{Z}_i, \mathbf{X}_i, U_i) L_b(V_j - V_i)(\mathbf{e}_j^{\mathrm{T}}, \mathbf{0}) + o_{\mathrm{P}}(c_n).$$

This completes the proof of Lemma A.4. □

LEMMA A.5. *Under Assumptions 1–5, we have* $\tilde{\mathbf{Z}}^{\mathrm{T}}(\mathbf{I}-\mathbf{S})(\mathbf{I}-\mathbf{S})^{\mathrm{T}}\tilde{\mathbf{Z}}/n \to \boldsymbol{\Sigma}$ *in probability and* $\hat{\boldsymbol{\Sigma}} = n(\tilde{\mathbf{Z}}^{\mathrm{T}}\tilde{\mathbf{Z}})^{-1}\tilde{\mathbf{Z}}^{\mathrm{T}}(\mathbf{I}-\mathbf{S})(\mathbf{I}-\mathbf{S})\tilde{\mathbf{Z}}(\tilde{\mathbf{Z}}^{\mathrm{T}}\tilde{\mathbf{Z}}^{\mathrm{T}})^{-1} \to \boldsymbol{\Sigma}$.

PROOF. The proof of the first result can be finished by arguments similar to those of Lemmas A.2–A.4, while the second one can be proved by arguments similar to Lemma 7.3 of Fan and Huang (2005). □

LEMMA A.6. *Under Assumptions 1–5, we have* $\tilde{\mathbf{Z}}^{\mathrm{T}}(\mathbf{I}-\mathbf{S})\mathbf{M}/n = O_{\mathrm{P}}(c_n^2)$.

PROOF. The proof follows (A.5) and an argument similar to that of Lemma 7.4 of Fan and Huang (2005). □

LEMMA A.7. $g(\cdot)$ *and* $h(\cdot)$ *are two continuous function vectors. Under Assumptions 1–5, we have* $\frac{1}{\sqrt{n}} \sum_{i=1}^{n} (\hat{\mathbf{Z}}_i - \mathbf{Z}_i) g(\mathbf{Z}_i) \varepsilon_i \to 0$ *and* $\frac{1}{\sqrt{n}} \sum_{i=1}^{n} (\hat{\mathbf{Z}}_i - \mathbf{Z}_i) \mathbf{X}_i^{\mathrm{T}} h(U_i) \varepsilon_i \to 0$ *in probability.*

PROOF. The proof follows from arguments similar to those of Lemma A.2. □

LEMMA A.8. *Under Assumptions 1–5, we have*

$$\tilde{\mathbf{Z}}^{\mathrm{T}}(\mathbf{I}-\mathbf{S})\boldsymbol{\varepsilon} = \sum_{i=1}^{n}\psi(\mathbf{Z}_i,\mathbf{X}_i,U_i)\mathbf{X}_i\{1+o_{\mathrm{P}}(1)\}\varepsilon_i + o(n^{1/2}),$$

where $\boldsymbol{\varepsilon} = (\varepsilon_1,\ldots,\varepsilon_n)^{\mathrm{T}}$.

PROOF. Note that $\tilde{\mathbf{Z}}^{\mathrm{T}}(\mathbf{I}-\mathbf{S})\boldsymbol{\varepsilon} = \sum_{i=1}^{n}\tilde{\mathbf{Z}}_i\{\varepsilon_i - (\mathbf{X}_i,0)(D_{u_i}W_{u_i}D_{u_i})^{-1}D_{u_i}\times W_{u_i}\boldsymbol{\varepsilon}\}$. By the same argument as those for (A.3), we have

$$n^{-1}D_u^{\mathrm{T}}W_u\boldsymbol{\varepsilon} = n^{-1}\sum_{i=1}^{n}\begin{pmatrix}\mathbf{X}_i \\ \dfrac{U_i-U}{h}\mathbf{X}_i\end{pmatrix}K_h(U_i-U)\varepsilon_i = f_u(U)E(\mathbf{X}|U)O_{\mathrm{P}}(c_n).$$

This formula along with (A.2) yields

$$(\mathbf{X}^{\mathrm{T}},0)(D_u^{\mathrm{T}}W_uD_u)^{-1}D_uW_u\boldsymbol{\varepsilon} = \mathbf{X}^{\mathrm{T}}\Gamma^{-1}(U)E(\mathbf{X}|U)O_{\mathrm{P}}(c_n).$$

A combination of these arguments with Lemma A.7 finishes the proof of Lemma A.8. □

PROOF OF THEOREM 1. Note that $\hat{\boldsymbol{\Theta}}_n$ can be expressed as $(\tilde{\mathbf{Z}}^{\mathrm{T}}\tilde{\mathbf{Z}})^{-1}\tilde{\mathbf{Z}}^{\mathrm{T}}(\mathbf{I}-\mathbf{S})\mathbf{Z}\boldsymbol{\Theta} + (\tilde{\mathbf{Z}}^{\mathrm{T}}\tilde{\mathbf{Z}})^{-1}\tilde{\mathbf{Z}}^{\mathrm{T}}(\mathbf{I}-\mathbf{S})\mathbf{M} + (\tilde{\mathbf{Z}}^{\mathrm{T}}\tilde{\mathbf{Z}})^{-1}\tilde{\mathbf{Z}}^{\mathrm{T}}(\mathbf{I}-\mathbf{S})\boldsymbol{\varepsilon}$. By Lemma A.8, the third term equals $\boldsymbol{\Sigma}^{-1}n^{-1}\sum_{i=1}^{n}\psi(\mathbf{Z}_i,\mathbf{X}_i,U_i)\varepsilon_i\{1+o_{\mathrm{P}}(1)\}+o_{\mathrm{P}}(n^{-1/2})$. The first term equals, via Lemma A.4,

$$\boldsymbol{\Sigma}^{-1}\left[\frac{\zeta_1 S_u^{-1}c_p b^{r+1}}{n(r+1)!}\sum_{i=1}^{n}\psi(\mathbf{Z}_i,\mathbf{X}_i,U_i)\{\boldsymbol{\xi}^{(r+1)}(V_i)\}^{\mathrm{T}}\boldsymbol{\beta}_0\right.$$
$$\left.+\frac{1}{n^2}\sum_{i=1}^{n}\sum_{j=1}^{n}\frac{1}{f_v(V_i)}\psi(\mathbf{Z}_i,\mathbf{X}_i,U_i)L_b(V_j-V_i)\mathbf{e}_j^{\mathrm{T}}\boldsymbol{\beta}_0\right].$$

By Lemma A.6 and (A.1), it follows that the second term of $\hat{\boldsymbol{\Theta}}_n$'s expression is of order $O(c_n^2)$ in probability. These arguments imply that

$$\hat{\boldsymbol{\Theta}}_n - \boldsymbol{\Theta}_0$$
$$= \boldsymbol{\Sigma}^{-1}\left[\frac{\zeta_1 S_u^{-1}c_p b^{r+1}}{n(r+1)!}\sum_{i=1}^{n}\psi(\mathbf{Z}_i,\mathbf{X}_i,U_i)\{\boldsymbol{\xi}^{(r+1)}(V_i)\}^{\mathrm{T}}\boldsymbol{\beta}_0\right.$$
$$+\frac{1}{n^2}\sum_{i=1}^{n}\sum_{j=1}^{n}\frac{1}{f_v(V_i)}\psi(\mathbf{Z}_i,\mathbf{X}_i,U_i)L_b(V_j-V_i)\mathbf{e}_j^{\mathrm{T}}\boldsymbol{\beta}_0$$
$$\left.+\frac{1}{n}\sum_{i=1}^{n}\psi(\mathbf{Z}_i,\mathbf{X}_i,U_i)\varepsilon_i\right]\{1+o_{\mathrm{P}}(1)\}.$$

This completes the proof of Theorem 1. □

PROOF OF THEOREM 3. By the definition of $\hat{\mathbf{\Psi}}(u)$, we have

$$\mathbf{H}\hat{\mathbf{\Psi}} = (D_u^{\mathrm{T}} W_u D_u)^{-1} D_u^{\mathrm{T}} W_u (Y - \hat{\mathbf{Z}}\hat{\mathbf{\Theta}}_n)$$
$$= I_1 + (D_u^{\mathrm{T}} W_u D_u)^{-1} D_u^{\mathrm{T}} W_u (\mathbf{Z} - \hat{\mathbf{Z}})\mathbf{\Theta}$$
$$+ (D_u^{\mathrm{T}} W_u D_u)^{-1} D_u^{\mathrm{T}} W_u \mathbf{Z} (\mathbf{\Theta} - \hat{\mathbf{\Theta}}_n) + R_n,$$

where $I_1 = (D_u^{\mathrm{T}} W_u D_u)^{-1} D_u W_u (Y - \mathbf{Z}\mathbf{\Theta})$ and $R_n = (D_u^{\mathrm{T}} W_u D_u)^{-1} D_u W_u (\mathbf{Z} - \hat{\mathbf{Z}})(\mathbf{\Theta} - \hat{\mathbf{\Theta}}_n)$. It is easy to show that $R_n = o(n^{-1/2})$ in probability. Note that

$$D_u^{\mathrm{T}} W_u (\mathbf{Z} - \hat{\mathbf{Z}})\mathbf{\Theta} = \left\{ \frac{\zeta_1 S_u^{-1} c_p b^{r+1}}{(r+1)!} \sum_{i=1}^n K_h(U_i - u) \right.$$
$$\times \begin{pmatrix} \mathbf{X}_i \{\boldsymbol{\xi}^{(r+1)}(V_i)\}^{\mathrm{T}} \boldsymbol{\beta} \\ h^{-1}(U_i - u)\mathbf{X}_i \{\boldsymbol{\xi}^{(r+1)}(V_i)\}^{\mathrm{T}} \boldsymbol{\beta} \end{pmatrix}$$
$$+ \frac{1}{n} \sum_{i=1}^n \sum_{j=1}^n f_v^{-1}(V_i) K_h(U_i - u) L_b(V_j - V_i)$$
$$\left. \times \begin{pmatrix} \mathbf{X}_i \mathbf{e}_j^{\mathrm{T}} \boldsymbol{\beta} \\ h^{-1}(U_i - u)\mathbf{X}_i \mathbf{e}_j^{\mathrm{T}} \boldsymbol{\beta} \end{pmatrix} \right\} \{1 + o_P(1)\}$$
$$\stackrel{\text{def}}{=} I_1' + I_2'.$$

It follows from (A.2) that

(A.8)
$$(D_u^{\mathrm{T}} W_u D_u)^{-1} I_1'$$
$$= \frac{\zeta_1 S_u^{-1} c_p b^{r+1}}{(r+1)!} \left\{ f_u(u) \Gamma^{-1}(u) \otimes \begin{pmatrix} 1 & \mu_1 \\ \mu_1 & \mu_2 \end{pmatrix} \right\}^{-1}$$
$$\times \sum_{i=1}^n \left\{ \begin{matrix} K_h(U_i - u)\mathbf{X}_i \{\boldsymbol{\xi}^{(r+1)}(V_i)\}^{\mathrm{T}} \boldsymbol{\beta} \\ K_h(U_i - u)\mathbf{X}_i h^{-1}(U_1 - u)\{\boldsymbol{\xi}^{(r+1)}(V_i)\}^{\mathrm{T}} \boldsymbol{\beta} \end{matrix} \right\} \{1 + o_P(1)\}$$

and

$$(D_u^{\mathrm{T}} W_u D_u)^{-1} I_2'$$
$$= \frac{1}{n^2} \sum_{j=1}^n \sum_{i=1}^n \left\{ f_u(u) \Gamma(u) \otimes \begin{pmatrix} 1 & \mu_1 \\ \mu_1 & \mu_2 \end{pmatrix} \right\}^{-1} K_h(U_i - u)$$
$$\times \left\{ \begin{matrix} E(\mathbf{X}_i | V = V_j)\mathbf{e}_j^{\mathrm{T}} \boldsymbol{\beta} \\ h^{-1}(U_i - u) E(\mathbf{X}_i | V = V_j)\mathbf{e}_j^{\mathrm{T}} \boldsymbol{\beta} \end{matrix} \right\} \{1 + o_P(1)\}$$

$$= \frac{\Gamma^{-1}(u)}{nf_u(u)(\mu_2 - \mu_1^2)} \sum_{i=1}^{n} K_h(U_i - u) E(\mathbf{X}_i | V = V_i) \mathbf{e}_i^{\mathrm{T}} \boldsymbol{\beta}$$

$$\otimes \begin{pmatrix} \mu_2 - \mu_1(U_i - u)/h \\ (U_i - u)/h - \mu_1 \end{pmatrix} \{1 + o_{\mathrm{P}}(1)\}.$$

Furthermore, (A.3) implies that

(A.9)
$$(D_u^{\mathrm{T}} W_u D_u)^{-1} D_u^{\mathrm{T}} W_u \mathbf{Z}(\boldsymbol{\Theta} - \hat{\boldsymbol{\Theta}}_n)$$
$$= \left\{ f_u(u) \Gamma(u) \otimes \begin{pmatrix} 1 & \mu_1 \\ \mu_1 & \mu_2 \end{pmatrix} \right\}^{-1}$$
$$\times \{n f_u(u) \Phi(u) \otimes (1, \mu_1)^{\mathrm{T}}\} (\boldsymbol{\Theta} - \hat{\boldsymbol{\Theta}}_n) \{1 + o_{\mathrm{P}}(1)\}$$
$$= \{\Gamma^{-1}(u) \Phi(u) \otimes (1, 0)^{\mathrm{T}}\} (\boldsymbol{\Theta} - \hat{\boldsymbol{\Theta}}_n) \{1 + o_{\mathrm{P}}(1)\}.$$

We therefore have $I_1 = (D_u^{\mathrm{T}} W_u D_u)^{-1} D_u^{\mathrm{T}} W_u \mathbf{M}_u + (D_u^{\mathrm{T}} W_u D_u)^{-1} D_u^{\mathrm{T}} W_u \boldsymbol{\varepsilon}$, where $\mathbf{M}_u = \boldsymbol{\alpha}(u)^{\mathrm{T}} \mathbf{X}$.

By the Taylor expansion and a direct simplification, we have

$$\mathbf{M} = \begin{pmatrix} \mathbf{X}_1^{\mathrm{T}} \boldsymbol{\alpha}(u) + (U_1 - u) \mathbf{X}_1^{\mathrm{T}} \boldsymbol{\alpha}'(u) + 2^{-1}(U_1 - u)^2 \mathbf{X}_1^{\mathrm{T}} \boldsymbol{\alpha}''(u) \\ \vdots \\ \mathbf{X}_n^{\mathrm{T}} \boldsymbol{\alpha}(u) + (U_n - u) \mathbf{X}_n^{\mathrm{T}} \boldsymbol{\alpha}'(u) + 2^{-1}(U_n - u)^2 \mathbf{X}_n^{\mathrm{T}} \boldsymbol{\alpha}''(u) \end{pmatrix} + o(h^2)$$

$$= D_u \begin{pmatrix} \boldsymbol{\alpha}(u) \\ h \boldsymbol{\alpha}'(u) \end{pmatrix} + \tfrac{1}{2} \begin{pmatrix} (U_1 - u)^2 \mathbf{X}_1^{\mathrm{T}} \boldsymbol{\alpha}''(u) \\ \vdots \\ (U_n - u)^2 \mathbf{X}_n^{\mathrm{T}} \boldsymbol{\alpha}''(u) \end{pmatrix} + o(h^2).$$

Hence,

(A.10)
$$I_1 = \left\{ \begin{pmatrix} \boldsymbol{\alpha}(u) \\ h\boldsymbol{\alpha}(u) \end{pmatrix} + \tfrac{1}{2} (D_u^{\mathrm{T}} W_u D_u)^{-1} D_u^{\mathrm{T}} W_u \right.$$
$$\left. \times \begin{pmatrix} (U_1 - u)^2 \mathbf{X}_1^{\mathrm{T}} \boldsymbol{\alpha}''(u) \\ \vdots \\ (U_n - u)^2 \mathbf{X}_n^{\mathrm{T}} \boldsymbol{\alpha}''(u) \end{pmatrix} + (D_u^{\mathrm{T}} W_u D_u)^{-1} D_u^{\mathrm{T}} W_u \boldsymbol{\varepsilon} \right\} \{1 + o_{\mathrm{P}}(h^2)\}.$$

It follows from (A.8)–(A.10) that $\sqrt{nh} \mathbf{H} \{\hat{\boldsymbol{\Psi}}(u_0) - \boldsymbol{\Psi}(u_0)\}$ can be represented as

$$\sqrt{nh} \left[ \frac{b^{r+1}}{n(r+1)!} \frac{\zeta_1^{\mathrm{T}} S_u^{-1} c_p}{(\mu_2 - \mu_1^2) f_u(u)} \Gamma^{-1}(U) \sum_{i=1}^{n} K_h(U_i - u) \right.$$
$$\otimes \begin{pmatrix} \{\mu_2 - \mu_1(U_i - u)/h\} \mathbf{X}_i \{\boldsymbol{\xi}^{(r+1)}(V_i)\}^{\mathrm{T}} \boldsymbol{\beta}_0 \\ \{(U_i - u)/h - \mu_1\} \mathbf{X}_i \{\boldsymbol{\xi}^{(r+1)}(V_i)\}^{\mathrm{T}} \boldsymbol{\beta}_0 \end{pmatrix}$$

$$+ \frac{\Gamma^{-1}(U)}{2(\mu_2 - \mu_1^2)f_u(u)} \sum_{i=1}^n K_h(U_i - u)(U_i - u)^2 \alpha''(u)$$

$$\otimes \begin{pmatrix} \{\mu_2 - \mu_1(U_i - u)/h\}\mathbf{X}_i\mathbf{X}_i^{\mathrm{T}} \\ \{(U_i - u)/h - \mu_1\}\mathbf{X}_i\mathbf{X}_i^{\mathrm{T}} \end{pmatrix} + o(h^2 + b^{r+1}) + O(n^{-1/2}) \Bigg]$$

$$+ \sqrt{nh}(D_u^{\mathrm{T}}W_u D_u)^{-1}D_u^{-1}W_u\boldsymbol{\varepsilon}$$

$$+ \frac{\sqrt{nh}\Gamma^{-1}(u)}{nf_u(u)(\mu_2 - \mu_1^2)} \sum_{i=1}^n K_h(U_i - u)$$

$$\times E(\mathbf{X}_i|V = V_i)\mathbf{e}_i^{\mathrm{T}}\boldsymbol{\beta} \otimes \begin{pmatrix} \mu_2 - \mu_1(U_i - u)/h \\ (U_i - u)/h - \mu_1 \end{pmatrix}$$

$$\times \{1 + o_{\mathrm{P}}(1)\}.$$

By an argument similar to that of Lemma A.8, we have

$$(D_u^{\mathrm{T}}W_u D_u)^{-1}D_u^{-1}W_u\boldsymbol{\varepsilon}$$
$$= \frac{\Gamma^{-1}(u)}{nf_u(u)(\mu_2 - \mu_1^2)} \sum_{i=1}^n K_h(U_i - u)\mathbf{X}_i\varepsilon_i$$
$$\otimes \begin{pmatrix} \mu_2 - \mu_1(U_i - u)/h \\ (U_i - u)/h - \mu_1 \end{pmatrix} \{1 + o_{\mathrm{P}}(1)\}.$$

The proof of Theorem 3 is completed. $\square$

PROOF OF THEOREM 5. The proof is similar to Theorems 3.1 and 3.2 of Fan and Huang (2005). We only give a sketch. We first prove that $n^{-1}RSS_1 = \sigma^2\{1 + o_{\mathrm{P}}(1)\}$.

By a procedure similar to that of Theorem 3.2 in Fan and Huang (2005), we can obtain that $n^{-1}RSS_{10} = n^{-1}\sum_{i=1}^n(Y_i - \hat{M}_{i0} - \hat{\boldsymbol{\Theta}}^{\mathrm{T}}\mathbf{Z}_i)^2 = \sigma^2\{1 + o_{\mathrm{P}}(1)\}$, where $\hat{M}_{i0}$ is the $i$th element of $\hat{\mathbf{M}}_0 = \mathbf{S}(\mathbf{Y} - \mathbf{Z}\hat{\boldsymbol{\Theta}})$. A direct calculation yields that

$$
\begin{aligned}
n^{-1}(RSS_1 - RSS_{10}) \\
= n^{-1}\sum_{i=1}^n \hat{\boldsymbol{\Theta}}^{\mathrm{T}}(\hat{\mathbf{Z}}_i - \mathbf{Z}_i)\{(Y_i - \hat{M}_i - \hat{\boldsymbol{\Theta}}^{\mathrm{T}}\hat{\mathbf{Z}}_i) + (Y_i - \hat{M}_{i0} - \hat{\boldsymbol{\Theta}}^{\mathrm{T}}\mathbf{Z}_i)\} \\
\text{(A.11)} \\
+ n^{-1}\sum_{i=1}^n (\hat{M}_i - \hat{M}_{i0})\{(Y_i - \hat{M}_i - \hat{\boldsymbol{\Theta}}^{\mathrm{T}}\hat{\mathbf{Z}}_i) \\
+ (Y_i - \hat{M}_{i0} - \hat{\boldsymbol{\Theta}}^{\mathrm{T}}\mathbf{Z}_i)\}.
\end{aligned}
$$

By (2.1), Theorem 2 and the Jensen inequality, we know that the first term in the right-hand side of (A.11) is bounded by

$$\max_{1 \leq i \leq n} \hat{\boldsymbol{\Theta}}^T |\hat{\mathbf{Z}}_i - \mathbf{Z}_i| \left[ \left\{ n^{-1} \sum_{i=1}^{n} (Y_i - \hat{M}_{i0} - \hat{\boldsymbol{\Theta}}^T \mathbf{Z}_i)^2 \right\}^{1/2} \right.$$
$$\left. + \max_{1 \leq i \leq n} \{|\hat{M}_i - \hat{M}_{i0}| + \hat{\boldsymbol{\Theta}}^T |\hat{\mathbf{Z}}_i - \mathbf{Z}_i|\} \right], \quad (A.12)$$

which is $o_P(1)$. A similar argument can show that the second term in the right-hand side of (A.11) is also $o_P(1)$. We therefore have $n^{-1}RSS_1 = \sigma^2\{1+o_P(1)\}$.

Furthermore, $RSS_0$ can be decomposed as $\{\mathbf{Y} - \hat{\mathbf{M}} - \hat{\mathbf{Z}}\hat{\boldsymbol{\Theta}} + \widetilde{\mathbf{Z}}(\hat{\boldsymbol{\Theta}} - \hat{\boldsymbol{\Theta}}_0)\}^T\{\mathbf{Y} - \hat{\mathbf{M}} - \hat{\mathbf{Z}}\hat{\boldsymbol{\Theta}} + \widetilde{\mathbf{Z}}(\hat{\boldsymbol{\Theta}} - \hat{\boldsymbol{\Theta}}_0)\} \stackrel{\text{def}}{=} RSS_1 + Q_1 + Q_2 + Q_3$, where $Q_1 = \{\widetilde{\mathbf{Z}}(\hat{\boldsymbol{\Theta}} - \hat{\boldsymbol{\Theta}}_0)\}^T\{\widetilde{\mathbf{Z}}(\hat{\boldsymbol{\Theta}} - \hat{\boldsymbol{\Theta}}_0)\}$, $Q_2 = (\mathbf{Y} - \hat{\mathbf{M}} - \hat{\mathbf{Z}}\hat{\boldsymbol{\Theta}})\{\widetilde{\mathbf{Z}}(\hat{\boldsymbol{\Theta}} - \hat{\boldsymbol{\Theta}}_0)\}$ and $Q_3 = \{\widetilde{\mathbf{Z}}(\hat{\boldsymbol{\Theta}} - \hat{\boldsymbol{\Theta}}_0)\}^T(\mathbf{Y} - \hat{\mathbf{M}} - \hat{\mathbf{Z}}\hat{\boldsymbol{\Theta}})$.

Recalling the expression of $\hat{\boldsymbol{\Theta}}_0$ and the result given in (A.1), we know that $n^{-1}\widetilde{\mathbf{Z}}^T\widetilde{\mathbf{Z}} \to \boldsymbol{\Sigma}$ in probability, and $Q_1 - n\boldsymbol{\Theta}^T\mathbf{A}^T\{\mathbf{A}\boldsymbol{\Sigma}^{-1}\mathbf{A}^T\}^{-1}\mathbf{A}\boldsymbol{\Theta} \to \sigma^2 \times \sum_{i=1}^{l} \omega_i \chi_{i1}^2$ in distribution. In an analogous way, we can show that $Q_2$ and $Q_3$ are asymptotic negligible in probability. These statements, along with the Slutsky theorem, imply that $2T_n - n\sigma^{-2}\boldsymbol{\Theta}^T\mathbf{A}^T\{\mathbf{A}\boldsymbol{\Sigma}^{-1}\mathbf{A}^T\}^{-1}\mathbf{A}\boldsymbol{\Theta} \to \sum_{i=1}^{l} \omega_i \chi_{i1}^2$ in distribution. Finally, following the lines of Rao and Scott (1981), we can prove that the distribution of $\varrho_n \sum_{i=1}^{l} \omega_i \chi_{i1}^2$ has the same approximate distribution as $\chi_l^2$, and complete the proof of Theorem 5. $\square$

**Acknowledgments.** The authors thank the Co-Editors, the former Co-Editors, an Associate Editor and two referees for their constructive comments, which substantially improved an earlier version of this paper.

Academy of Mathematics and Systems Science  
Chinese Academy of Sciences  
Beijing 100080  
China  
and  
Department of Statistics  
Shanghai University of Finance and Economics  
Shanghai 200433  
China  
E-mail: yzhou@amss.ac.cn

Department of Biostatistics  
and Computational Biology  
University of Rochester  
Rochester, New York 14642  
USA  
E-mail: hliang@bst.rochester.edu